\documentclass[11pt]{article}
\usepackage{graphicx,epsfig}
\usepackage{amsmath,amssymb,euscript}
\usepackage{latexsym}
\usepackage{color}
\newtheorem{theorem}{\bf Theorem}[section]

\newtheorem{lemma}[theorem]{\bf Lemma}

\newtheorem{corollary}[theorem]{\bf Corollary}
\newtheorem{remark}[theorem]{\bf Remark}

\newcommand{\Proofend}{\hfill$\diamondsuit$}

\def\NN{{\mathbb N}}
\def\RR{{\mathbb R}}

\def\PP{{\mathbb P}}

\def\C{{\mathcal{C}}}

\def\O{{\mathcal{O}}}

\def\tlambda{{\tilde{\lambda}}}

\graphicspath{ {./} }

\begin{document}
\title{A generalization of Floater--Hormann interpolants}
\author{Woula Themistoclakis\thanks{C.N.R. National
        Research Council of Italy,
        IAC Institute for Applied Computing ``Mauro Picone'',  via P. Castellino, 111, 80131 Napoli, Italy.
        woula.themistoclakis@cnr.it.
This author was partially supported by GNCS-INDAM and the research has been accomplished within RITA - Research ITalian network on Approximation - and UMI - T.A.A. working group }
        \and Marc Van Barel\thanks{KU Leuven, Department of Computer Science, KU Leuven,
Celestijnenlaan 200A,
B-3001 Leuven (Heverlee), Belgium. marc.vanbarel@cs.kuleuven.be.
This author was supported by the Research Council KU Leuven, C1-project C14/17/073 and by the Fund for Scientific Research–Flanders (Belgium), EOS Project no 30468160 and project G0B0123N.}}
\maketitle

\begin{abstract}
In this paper the interpolating rational functions introduced by Floater and Hormann are generalized leading to a whole new family of rational functions depending on  $\gamma$, an additional positive integer parameter.  For $\gamma = 1$, the original Floater--Hormann interpolants are obtained. When $\gamma>1$ we prove  that the new rational functions share a lot of the nice properties of the original Floater--Hormann functions. Indeed, for any configuration of nodes in a compact interval, they have no real poles, interpolate the given data,  preserve the polynomials
up to a certain fixed degree, and have a barycentric-type representation.
Moreover, we estimate the associated Lebesgue constants in terms of the minimum ($h^*$) and maximum ($h$) distance between two consecutive nodes. It turns out that,
in contrast to the original Floater-Hormann interpolants, for all $\gamma > 1$ we get
 uniformly bounded Lebesgue constants in the case of equidistant and quasi-equidistant nodes configurations (i.e., when $h\sim h^*$). For such configurations, as the number of nodes tends to infinity, we prove that the new interpolants ($\gamma>1$) uniformly converge to the interpolated function $f$, for any continuous function $f$ and all $\gamma>1$. The same is not ensured by the original FH interpolants ($\gamma=1$). Moreover, we provide uniform and pointwise estimates of the approximation error for functions having different degrees of smoothness.
  Numerical experiments illustrate the theoretical results and show a better error profile
for less smooth functions compared to the original Floater-Hormann interpolants.
\end{abstract}

\section{Introduction}

In this paper, we consider the problem of interpolating a function $f(x)$ on a finite interval $[a,b]$, given its values $f(x_i)$ in
$n+1$  nodes $a = x_0 < x_1< \cdots < x_{n-1} < x_n=b$.
Without loss of generality the interval $[-1,+1]$ can be taken.
If one can choose the position of the nodes $x_i$ in the interval, an analytic function can be approximated by polynomials interpolating
at the Chebyshev points of the first or second kind, leading to exponential convergence.
The speed of convergence is determined by the largest Bernstein ellipse that can be taken in the analytic domain of the function.
For differentiable functions, the convergence is algebraic where the speed of convergence is determined by the smoothness
of the function. For further details we refer the reader to the book of Trefethen \cite{b495}.

If the nodes can not be freely chosen, the problem becomes much harder. E.g., when the nodes are equidistant in the interval
$[-1,+1]$ and we want to approximate the function $1/(1+25x^2)$, the Runge phenomenon occurs and the approximation
error becomes very large in the neighborhood of the endpoints of the interval when the number of nodes increases.
In  \cite{huybrechs2023aaa} Huybrechs and Trefethen compare several methods to approximate a function when the nodes are equidistant.
One of these methods uses the Floater--Hormann (briefly FH) interpolating rational functions \cite{floater2007barycentric}.
This method is valid for any configuration of the nodes but turns out to be very useful for equidistant configurations.
Generalizing \cite{berrut1988rational},  Floater and Hormann introduced a blended form of interpolating polynomials of fixed degree $0< d\le n$, leading to the rational function  \cite{floater2007barycentric}
\begin{equation}\label{eq:001}
r(x) = \frac{\sum_{i=0}^{n-d} \lambda_i(x) p_i(x)}{\sum_{i=0}^{n-d} \lambda_i(x)}
\end{equation}
where, for all $i=0,\ldots, (n-d)$,
\begin{equation}\label{eq:002}
\lambda_i(x) = \frac{(-1)^i}{(x-x_i)(x-x_{i+1})\ldots (x-x_{i+d})},
\end{equation}
and $p_i(x)$ denotes the unique polynomial of degree $\leq d$ interpolating $f$ at the $(d+1)$ points $x_i<x_{i+1}<\ldots<x_{i+d}$.

Note that in the limit case $d=0$, we get $p_i(x)=f(x_i)$, $\lambda_i(x)=(-1)^i(x-x_i)^{-1}$, and (\ref{eq:001}) yields the Berrut rational interpolants studied in \cite{berrut1988rational}. By taking any $1\le d\le n$,
Floater and Hormann proved that the approximant $r(x)$ has no poles on the real line, it coincides with $f$ on the set of nodes
$X_n=\{x_0,x_1,\ldots,x_n\}$, and, at any $x\not\in X_n$, it admits a barycentric representation allowing efficient and stable computations \cite{camargo2016numerical}.
Moreover, for any fixed $d\in\NN$, as $h = \max_{1\le i\le n} (x_{i+1}-x_i)\to 0$ (and hence $n\to\infty$), regardless of the distribution of nodes, the FH approximation error behaves as follows \cite[Thm. 2]{floater2007barycentric}
\begin{equation}\label{err-FH}
\| r-f \|_\infty = \O(h^{d+1}), \qquad \forall f \in \C^{d+2}([a,b]),
\end{equation}
where, as usual, $\|f\|_\infty=\max_{x\in [a,b]}|f(x)|$.

Therefore the FH interpolants, which for $d=0$ reduce to Berrut interpolants \cite{floater2007barycentric}, are able to produce arbitrarily high approximation orders provided that the parameter $d$ is large enough. However, in the important case of equidistant or quasi--equidistant nodes, it has been proved that the Lebesgue constants of FH interpolants grow logarithmically with $n$, but exponentially with $d$ \cite{bos2012lebesgue, hormann2012barycentric}. Hence, increasing $d$ too much is not always advisable.

For an overview of linear barycentric rational interpolation, we refer the interested reader to the paper of Berrut and Klein
\cite{berrut2014recent}. This overview also describes a generalization of the FH interpolant developed
by Klein \cite{klein2013extension} in the case of equidistant nodes.
See also \cite{camargo2017comparison,camargo2017stability}.

In this paper we are also going to generalize the method of Floater and Hormann but in a different way.
For any distribution of nodes, we define a whole family of linear rational approximants that are denoted by $\tilde r(x)$ and depend, besides $d$, on an additional parameter $\gamma \in \NN = \{1,2,3,\ldots\}$. When $\gamma = 1$, $\tilde r(x)$ reduces to the original FH interpolant $r(x)$.
When $\gamma=2$ and $d=0$, $\tilde r(x)$ reduces to the approximants studied in \cite{zhang2022rational}. In this paper, for brevity,  we only consider the case of arbitrary $1\le d\le n$, the case $d=0$ will be investigated in future work.

Similarly to the original FH interpolants $r(x)$, we show that, for all $\gamma>1$, also $\tilde r(x)$ has no real poles, interpolates the data, and has a barycentric type
representation. However, the main results of the paper concern the case of equidistant or quasi-equidistant configurations of nodes, where we find some novelty by taking $\gamma>1$. First of all, we prove that
\begin{equation}\label{limit}
\lim_{n\to\infty}\|\tilde r -f\|_\infty=0,\qquad \forall f\in C([a,b]),\qquad \forall d\in\NN,\qquad \forall \gamma>1,
\end{equation}
which is not ensured when $\gamma=1$. Moreover, in contrast to the original FH interpolants, we prove that for equidistant or quasi-equidistant nodes, the Lebesgue constants corresponding to any $1\le d\le n$ and $\gamma>1$ are uniformly bounded both in $n$ and $\gamma$, but still grow exponentially with $d$.

With respect to  the approximation rate, for equidistant or quasi-equidistant nodes, we show that
\begin{equation}\label{err-FH-new}
\| \tilde r-f \|_\infty = \O(h^{s}), \qquad \forall f \in \C^{s}([a,b]),\qquad 1\le s\le d+1,
\end{equation}
holds as $n\to\infty$, for arbitrarily fixed parameters $d\in\NN$ and $\gamma>s+1$. Hence, generalized FH interpolants also provide arbitrarily high convergence rates. Moreover, making the comparison with (\ref{err-FH}), if we take $s=d+1$ in (\ref{err-FH-new}) then we  get that the generalized FH functions $\tilde r(x)$, with parameter $\gamma>d+2$, also can reach the convergence order $\O(h^{d+1})$ but supposing that $f\in C^{d+1}([a,b])$  instead of $f\in C^{d+2}([a,b])$ .

In addition, we estimate the error of generalized FH interpolants also for functions with a non-integer smoothness degree. More precisely, in the class $Lip_\alpha ([a,b])$ of H\"older continuous functions with exponent $0<\alpha\le 1$, we prove that
\begin{equation}\label{err-FH-new1}
\| \tilde r-f \|_\infty = \O(h^{\alpha}), \qquad \forall f \in Lip_{\alpha}([a,b]),\qquad \forall\gamma>\alpha +1.
\end{equation}
Moreover, in  the class $C^{s,\alpha}([a,b])$  of functions that are $s$-times continuously differentiable and have the $s$--th derivative H\"older continuous with exponent $0<\alpha\le 1$, the generalized FH function $\tilde r(x)$ corresponding to any $d\in\NN$ satisfies
\begin{equation}
\label{err-FH-new2}
\| \tilde r-f \|_\infty = \O(h^{s+\alpha}),\quad \qquad \forall f \in C^{s,\alpha}([a,b]),\qquad\quad  \forall \gamma>s+\alpha+1.\\
\end{equation}
 for all integers $1\le s\le d$. 

Finally, we consider the pointwise error and show that, even in the case of almost everywhere continuous functions $f$ with some isolated discontinuities, the approximation orders displayed in  (\ref{err-FH-new})--(\ref{err-FH-new2}) are locally preserved and continue to hold in all compact subintervals $I\subset [a,b]$ where  $f$ has the described smoothness degree (i.e., $f \in C^{s}(I)$, $f \in Lip_{\alpha}(I)$, and $f \in C^{s,\alpha}(I)$, resp.)

The numerical experiments confirm the theoretical estimates and show that less restrictive assumptions on $\gamma$ could be possible for getting the previous error trend. Moreover, in comparison with the original FH interpolants,  the generalized ones exhibit a much better error profile when the interpolated function is less smooth.

The paper is organized as follows. In Section \ref{sec:def} the generalized FH rational interpolants are presented.
In Section \ref{sec:prop} we state and prove several properties of these new interpolants which are similar to those of the original
FH interpolants.
In Section \ref{sec:LC} we focus on the case of equidistant or quasi-equidistant nodes and estimate the behaviour of the associated Lebesgue constants. This section concludes with two subsections: in the former the proof and the necessary technical lemmas are given, in the latter a useful related result is stated. In Section \ref{sec:er} the convergence theorem and all the error estimates are given.
In Section \ref{sec:ex} we illustrate several numerical examples comparing generalized and original FH approximants. Finally,
Section \ref{sec:con}  gives the conclusion of our paper.

\section{Generalized Floater--Hormann interpolants}\label{sec:def}
Let
\begin{equation}\label{xk}
a=x_0<x_1<\ldots < x_{n-1}<x_n=b,\qquad\qquad  n\in\NN,
\end{equation}
be any sequence of nodes where we assume the function $f:[a,b]\to \RR$ has been sampled.
The generalized FH approximation of $f$ is defined very similarly to (\ref{eq:001}) and (\ref{eq:002}).
For any integer $0\le d\le n$, it is also a blended form of the polynomial interpolants $p_i$ of degree at most $d$.
However, the blending functions depend on an additional parameter $\gamma \in \NN$ and are defined as follows
\begin{equation}\label{eq:004}
\tlambda_i(x) = \frac{(-1)^{i \gamma}}{(x-x_i)^\gamma(x-x_{i+1})^\gamma\ldots (x-x_{i+d})^\gamma},\qquad i=0,\ldots, n-d.
\end{equation}
Hence, for arbitrarily fixed $\gamma\in \NN$ and $d\in\{0,\ldots, n\}$, the generalized FH approximation of $f$ is given by
\begin{equation}\label{eq:003}
\tilde{r}(x) = \frac{\sum_{i=0}^{n-d} \tlambda_i(x) p_i(x)}{\sum_{i=0}^{n-d} \tlambda_i(x)}, \qquad x\in\RR
\end{equation}
where $\tlambda_i(x)$ is defined in (\ref{eq:004}) and $p_i(x)$ is the polynomial of degree $\leq d$ interpolating $f$ at the $(d+1)$ nodes
$x_i<x_{i+1}<\ldots<x_{i+d}$.

We point out that the function $\tilde r(x)$ depends on $f$, on $n$, on the nodes (\ref{xk}), and on two integer parameters: $0\le d\le n$ and $\gamma\ge 1$. Sometimes we also use the notation $\tilde r(f,x)=\tilde r(x)$ and $\tilde r_d(f,x)=\tilde r(x)$ in order to highlight the dependence on $f$ and $d$.

Note that the original FH interpolants are a special case of the generalized ones with parameter $\gamma = 1$ (compare (\ref{eq:002}) and (\ref{eq:004})).

If we multiply both the numerator and denominator of (\ref{eq:003}) by the polynomial
\begin{equation}\label{pol-pi}
\pi(x)=(-1)^{\gamma (n-d)}\prod_{k=0}^n(x-x_k)^\gamma,
\end{equation}
then we get
\begin{equation}\label{r-pol}
\tilde r(x)=\frac{\sum_{i=0}^{n-d}\tilde \mu_i(x)p_i(x)}{\sum_{i=0}^{n-d}\tilde \mu_i(x)}
\end{equation}
where we set
\begin{equation}\label{mu-tilde}
\tilde \mu_i(x)=\tlambda_i(x)\pi(x)=\prod_{k=0}^{i-1}(x-x_k)^\gamma \prod_{k=i+d+1}^n (x_k-x)^\gamma
\end{equation}
being understood that $\prod_{k=n_1}^{n_2}a_k=1$ whenever the product is empty, i.e., if $n_1>n_2$.

Equation (\ref{r-pol}) yields the generalized FH approximant $\tilde r(x)$ as a quotient of two polynomials
\begin{equation}
\label{PQ}
P(x)=\sum_{i=0}^{n-d}\tilde \mu_i(x)p_i(x)\\
\qquad \quad\mbox{and}\qquad\quad
Q(x)=\sum_{i=0}^{n-d}\tilde \mu_i(x)
\end{equation}
where the maximum degree is  $\gamma (n-d)+d$ at the numerator and $\gamma(n-d)$ at the denominator, i.e., for all $\gamma\in\NN$, $\tilde r(x)$ is a rational function of type $(\gamma n- (\gamma -1)d, \ \gamma n-\gamma d)$.

\section{Properties}\label{sec:prop}
In this section, we are going to prove that, for any choice of the integer parameter $\gamma>1$, the generalized FH function $\tilde r(x)$  has similar properties to the original $r(x)$.
\subsection{No poles on the real line}
In the case $\gamma=1$ the polynomials $\tilde\mu_i(x)$ defined in (\ref{mu-tilde}) were investigated by Floater and Hormann in \cite[Thm. 1]{floater2007barycentric}. Using their result, the following lemma can be easily deduced:
\begin{lemma}\label{lem}
Let $n\in\NN$, $d\in\{0,1,\ldots, n\}$, and $i\in\{0,1,\ldots,n-d\}$ be arbitrarily fixed.

For any $\gamma\in\NN$ that is even, and $\forall x\in\RR$,  we have
\begin{equation}\label{mu-even}
\tilde\mu_i(x)\left\{\begin{array}{cccc}
& =0 & &\mbox{if} \quad x\in \left\{x_k \ :\ \ 0\le k< i\ \bigvee \ i+d<k\le n\right\}\\
& >0 & & \mbox{otherwise}
\end{array}\right.
\end{equation}
If $\gamma\in\NN$ is odd, we distinguish the following cases:
\begin{itemize}
\item In the case $x\in\{x_0,\ldots, x_n\}$ we have
\begin{equation}\label{mu-xk}
\tilde\mu_i(x_k)\left\{\begin{array}{cccl}
& >0 & &\mbox{if} \quad (k-d)\le i\le k\\
& =0 & & \mbox{otherwise}
\end{array}\right.
\end{equation}
\item In the case $x\in [a,b]-\{x_0,\ldots, x_n\}$, if $x_\ell<x<x_{\ell +1}$ for any $\ell=0,\ldots, n-1$, then we have the following
\begin{itemize}
\item[(A)] If $(\ell-d+1)\le i\le\ell$ then $\tilde\mu_i(x)>0$
\item[(B)] Suppose $\ell-d\ge 0$, for $i=0,1,\ldots,(\ell-d)$,  the sequence $\tilde\mu_i(x)$ has alternate signs, ending with $\tilde\mu_{\ell-d}(x)>0$, and it has increasing absolute values, i.e.,
    \begin{equation}\label{case-less}
\tilde\mu_{\ell-d-k}(x)> -\tilde\mu_{\ell-d-k-1}(x)>0,\qquad k=0,2,4,\ldots
    \end{equation}
\item[(C)] For $i=(\ell+1),\ldots, n$,  the sequence $\tilde\mu_i(x)$ has alternate signs, starting with $\tilde\mu_{\ell+1}(x)>0$, and it has decreasing absolute values, i.e.,
    \begin{equation}\label{case-more}
\tilde\mu_{\ell+k}(x)> -\tilde\mu_{\ell+k+1}(x)>0,\qquad k=1,3,5,\ldots
    \end{equation}
\end{itemize}
\item In the case $x<a$,  the sequence $\tilde\mu_i(x)$, for $i=0,\ldots, (n-d)$, has alternate sign, starting with $\tilde\mu_0(x)>0$, and it has decreasing absolute values, i.e.,
    \begin{equation}\label{case-a}
\tilde\mu_{k}(x)> -\tilde\mu_{k+1}(x)>0,\qquad k=0,2,4,\ldots
    \end{equation}
    \item In the case $x>b$,  the sequence $\tilde\mu_i(x)$, for $i=0,\ldots, (n-d)$, has alternate sign, ending with $\tilde\mu_{n-d}(x)>0$, and it has increasing absolute values, i.e.,
    \begin{equation}\label{case-b}
\tilde\mu_{n-d-k}(x)> -\tilde\mu_{n-d-k-1}(x)>0,\qquad k=0,2,4,\ldots
    \end{equation}
\end{itemize}
\end{lemma}
From the previous lemma we deduce the following result that generalizes \cite[Thm. 1]{floater2007barycentric}.
\begin{theorem}\label{th-poles}
For all integers $n,\gamma\in\NN$ and $0\le d\le n$, the generalized FH rational function $\tilde r(x)$ has no real poles.
\end{theorem}
{\it Proof of Theorem \ref{th-poles}}\newline
Recalling (\ref{r-pol}), it is sufficient to prove that
\begin{equation}\label{tesi}
Q(x)=\sum_{i=0}^{n-d}\tilde \mu_i(x)>0, \qquad \forall x\in\RR.
\end{equation}
This follows from (\ref{mu-even}) in the case that $\gamma\in\NN$ is even. If $\gamma\in\NN$ is odd, by (\ref{mu-xk}) we get
\[
Q(x_k)=\sum_{i=max\{0, (k-d)\}}^k \tilde\mu_i(x_k)>0, \qquad k=0,1,\ldots,n,
\]
while (\ref{case-a}) and (\ref{case-b}) imply
\[
Q(x)=\left\{
\begin{array}{lr}
\displaystyle\sum_{\scriptsize {\begin{array}{c}{k=0,2,4,...}\\ [-.1cm]
 k< n-d\end{array}}}[\tilde\mu_{k}(x)+\tilde \mu_{k+1}(x)]>0,
 & \qquad\forall x<a\\ [1cm]
\displaystyle\sum_{\scriptsize {\begin{array}{c}{k=0,2,4,...}\\ [-.1cm]
 k< n-d\end{array}}}[\tilde\mu_{n-d-k}(x)+\tilde \mu_{n-d-k-1}(x)]>0,
 & \qquad\forall x>b,
\end{array}\right.
\]
in the case that $(n-d)$ is odd, and
\[
Q(x)=\left\{
\begin{array}{ll}
\displaystyle\sum_{\scriptsize {\begin{array}{c}{k=0,2,4,...}\\ [-.1cm]
 k< n-d\end{array}}}[\tilde\mu_{k}(x)+\tilde \mu_{k+1}(x)]+ \tilde\mu_{n-d}(x)>0, &\qquad \forall x<a\\ [1cm]
\displaystyle\sum_{\scriptsize {\begin{array}{c}{k=0,2,4,...}\\ [-.1cm]
 k< n-d\end{array}}}[\tilde\mu_{n-d-k}(x)+\tilde \mu_{n-d-k-1}(x)]+\tilde\mu_0(x)>0, &\qquad \forall x>b,
\end{array}\right.
\]
when $(n-d)$ is even.
Hence, it remains to prove (\ref{tesi}) only in the case that $\gamma\in\NN$ is odd and $x_\ell<x<x_{\ell+1}$ for some $\ell=0,1,\ldots, n-1$. In such a case, we write
\begin{eqnarray*}
Q(x)&=&\sum_{0\le i\le\ell-d}\tilde\mu_i(x)+\sum_{\ell-d+1\le i\le\ell}\tilde\mu_i(x)+\sum_{\ell+1\le i\le n}\tilde\mu_i(x)\\ [.1in]
&=:&Q_1(x)+Q_2(x)+Q_3(x)
\end{eqnarray*}
being understood that $\sum_{n_1\le i\le n_2}a_i=0$ in case of empty summation, i.e., if $n_1>n_2$.

Finally, we observe that whenever the previous summation $Q_i(x)$ are non empty, they are positive by virtue of Lemma \ref{lem} (cf. (A)--(C)). Since at least one term of the summation $Q_i(x)$ is not empty, we have proven the theorem.
\Proofend
\begin{remark}
We remark that the previous proof also states that for all $n,\gamma\in\NN$, and $1\le d\le n$,  if $x\in ]x_\ell, \ x_{\ell+1}[$, with $\ell =0,\ldots, n-1$, then we have
\begin{equation}\label{eq-rem}
|Q(x)|=\left|\sum_{i=0}^{n-d}\tilde\mu_i(x)\right|\ge
\sum_{i\in I_\ell}\tilde\mu_i(x)>\tilde\mu_j(x)>0,\qquad \forall j\in I_\ell
\end{equation}
where
\begin{equation}\label{Il}
I_\ell=\left\{i\in\{0,\ldots, (n-d)\} \ : \ \ell-d+1\le i\le \ell\right\}
\end{equation}
Note that $I_\ell$ is a non empty set  since $d\ge1$.
Moreover, if we divide all terms in (\ref{eq-rem}) by $|\pi(x)|$, with $\pi(x)$ defined in (\ref{pol-pi}), then, by virtue of (\ref{mu-tilde}), we obtain
\begin{equation}\label{rem-lambda}
\left|\sum_{i=0}^{n-d}\tlambda_i(x)\right|\ge
\sum_{i\in I_\ell}\tlambda_i(x)>\tlambda_j(x)>0,\qquad \forall j\in I_\ell,\qquad \forall x\in ]x_\ell,\ x_{\ell +1}[
\end{equation}
\end{remark}
\subsection{Interpolation of the data}
In the case $\gamma=1$ it is known that the FH rational function $r(x)$ is equal to $f(x)$ if $x$ is one of the nodes (\ref{xk}). Such interpolation property remains valid for the generalized FH approximation.
\begin{theorem}\label{th-interp}
For all $n,\gamma\in\NN$ and any integer $0\le d\le n$, we have
\begin{equation}
\label{eq-interp}
\tilde r(f,x_k)=f(x_k), \qquad k=0,1,\ldots, n
\end{equation}
\end{theorem}
{\it Proof of Theorem \ref{th-interp}}
Set
\begin{equation}\label{Jk}
J_k=\left\{i\in\{0,\ldots, (n-d)\} \ : \ k-d\le i\le k\right\},
\qquad k=0,\ldots,n
\end{equation}
we note that
\[
p_i(x_k)=f(x_k),\qquad \forall i\in J_k.
\]
Moreover, by Lemma \ref{lem} we have
\[\tilde\mu_i(x_k)=0, \qquad \forall i\not\in J_k
\]
Consequently, we get
\[
\tilde r(f,x_k)=\frac{\sum_{i=0}^{n-d}\tilde \mu_i(x_k)p_i(x_k)}{\sum_{i=0}^{n-d}\tilde \mu_i(x_k)}=
\frac{\sum_{i\in J_k}\tilde \mu_i(x_k)p_i(x_k)}{\sum_{i\in J_k}\tilde \mu_i(x_k)}=f(x_k).
\]
\Proofend
\subsection{Preservation of polynomials}
Similarly to the classical FH interpolants, also the generalized ones reduce to the identity on the set $\PP_d$ of all polynomials of degree at most $d$.

Indeed, for all $f\in\PP_d$ we have
\[
p_i(x)=f(x), \qquad i=0,\ldots, n-d, \qquad \forall x\in [a,b].
\]
Consequently, by (\ref{r-pol}), we deduce that
\begin{equation}\label{eq-inva}
\tilde r(f,x)=f(x), \qquad \forall x\in [a,b], \qquad \forall f\in\PP_d,
\end{equation}
holds for any $n,\gamma\in\NN$, and $d\in\{0,\ldots, n\}$.
\subsection{Barycentric-type representation}
We recall that the classical FH interpolants are rational functions of type $(n,\ n-d)$, hence they can be expressed in barycentric form.
In this subsection, we give a barycentric-type expression for the generalized FH interpolants defined by
(\ref{eq:003}) and (\ref{eq:004}), for any parameter $\gamma\in\NN$.

The Lagrange representation for the interpolating polynomial $p_i$ is
$$
p_i(x) = \sum_{k=i}^{i+d} f(x_k)\prod_{s=i,s\neq k}^{i+d} \frac{x-x_s}{x_k-x_s} .
$$
Combining this with (\ref{eq:004}), we get
\begin{eqnarray*}
\tlambda_i(x) p_i(x) &=& (-1)^{i \gamma}\prod_{s=i}^{i+d} \frac{1}{(x-x_s)^\gamma} \left[\sum_{k=i}^{i+d} f(x_k)
\prod_{s=i,s\neq k}^{i+d}
\frac{x-x_s}{x_k-x_s}\right] \\
&=& (-1)^{i \gamma}\sum_{k=i}^{i+d} \frac{f(x_k)}{(x-x_k)^\gamma} \prod_{s=i,s\neq k}^{i+d} \frac{1}{(x_k-x_s)(x-x_s)^{\gamma-1}}.
\end{eqnarray*}
Hence, we obtain
$$
\sum_{i=0}^{n-d} \tlambda_i(x) p_i(x) = \sum_{i=0}^{n-d} (-1)^{i \gamma}\sum_{k=i}^{i+d} \frac{f(x_k)}{(x-x_k)^\gamma} \prod_{s=i,s\neq k}^{i+d} \frac{1}{(x_k-x_s)(x-x_s)^{\gamma -1}},
$$
and changing the order of the summations, we get
$$
\sum_{i=0}^{n-d} \tlambda_i(x) p_i(x) = \sum_{k=0}^n \frac{f(x_k)}{(x-x_k)^\gamma} \left[\sum_{i \in J_k} (-1)^{i\gamma} \prod_{j=i,j\neq k}^{i+d} \frac{1}{(x_k-x_j)(x-x_j)^{\gamma -1}}\right]
$$
where we recall $J_k = \left\{  i \in \{0,1,\ldots,n-d\} : k-d \leq i \leq k \right\}$ has been introduced in (\ref{Jk}).

Hence, defining
\begin{equation}\label{wk}
w_k(x) = \sum_{i \in J_k} (-1)^{i\gamma} \prod_{s=i,s\neq k}^{i+d} \frac{1}{(x_k-x_s)(x-x_s)^{\gamma -1}},\qquad \gamma\in\NN
\end{equation}
we can write
\begin{equation}\label{bar-num}
\sum_{i=0}^{n-d} \tlambda_i(x) p_i(x) = \sum_{k=0}^n \frac{f(x_k)}{(x-x_k)^\gamma}w_k(x).
\end{equation}
Similarly, we obtain
\begin{equation}\label{bar-den}
\sum_{i=0}^{n-d} \tlambda_i(x) = \sum_{k=0}^n \frac{1}{(x-x_k)^\gamma} w_k(x),
\end{equation}
that is (\ref{bar-num}) in the case of the unit function $f(x)=1$, $x\in [a,b]$.

In conclusion, the generalized FH interpolant (\ref{eq:003}) can be expressed in the following barycentric-type form
\begin{equation}\label{eq-bar}
\tilde{r}(x) = \frac{\sum_{k=0}^n \frac{w_k(x)}{(x-x_k)^\gamma} \ f(x_k)}{\sum_{k=0}^n \frac{ w_k(x)}{(x-x_k)^\gamma}},\qquad \gamma\in\NN
\end{equation}
Note that in the case $\gamma=1$, this yields the classical barycentric form
\begin{equation}\label{bar-classic}
{r}(x) = \frac{\sum_{k=0}^n \frac{ w_k}{(x-x_k)}f(x_k)}{\sum_{k=0}^n \frac{w_k}{(x-x_k)}},\qquad\quad
w_k= \sum_{i \in J_k} (-1)^{i} \prod_{s=i,s\neq k}^{i+d} \frac{1}{(x_k-x_s)},
\end{equation}
where the weights $w_k$ can be computed in advance and where, in the denominators,
we have a factor $(x-x_k)$ instead of $(x-x_k)^\gamma$.

Note that for $\gamma > 1$ the weights can not be precomputed.
Hence,  to evaluate the new interpolant, in $m$ $x$-values $\O(m n d^2)$ FLOPS are needed while for the classical barycentric form
the weights can be computed beforehand using (\ref{bar-classic}) in $\O(n d^2)$ FLOPS.
Using a more complicated pyramid algorithm \cite{hormann2016pyramid}, this can even be reduced to $\O(n d)$ FLOPS.
 Evaluating the classical barycentric form in $m$ $x$-values costs an additional
$\O(mn)$ FLOPS.
Hence, it is more efficient to evaluate the classical FH interpolant  in comparison to the new one.
However, the new approximant exhibits better performance with respect to the error, especially for less smooth functions, as
will be shown in the numerical examples.
\section{Lebesgue constants}\label{sec:LC}
Set for brevity
\begin{equation}\label{D}
D(x)=\sum_{i=0}^{n-d}\tlambda_i(x).
\end{equation}
By (\ref{eq:003}) and (\ref{bar-num}) we get
\[
\tilde r(f, x)=\sum_{k=0}^n f(x_k)\frac{w_k(x)}{(x-x_k)^\gamma D(x)}
\]
i.e., defining
\begin{equation}\label{bk}
b_k(x)=\left\{\begin{array}{cll}
1 & \mbox{if}& x=x_k\\
0 & \mbox{if}& x\in \{x_0,\ldots, x_n\}-\{x_k\}\\
\displaystyle \frac{w_k(x)}{(x-x_k)^\gamma D(x)} & \mbox{if}& x\not\in\{x_0,\ldots, x_n\}
\end{array}\right.\qquad k=0,\ldots, n,
\end{equation}
we can write
\begin{equation}\label{r-bk}
\tilde r(f,x)=\sum_{k=0}^n f(x_k) b_k(x), \qquad x\in\RR .
\end{equation}
The Lebesgue constant and function of the generalized FH interpolants at the nodes (\ref{xk}) are given by
\[
\Lambda_n=\sup_{x\in[a,b]}|\Lambda_n(x)|,\qquad \quad \Lambda_n(x)=\sum_{k=0}^n |b_k(x)|.
\]
Their behaviour as $n\rightarrow\infty$ is an important measure for the conditioning of the problem, being well--known that
\[
|\tilde r(f,x)-\tilde r(F,x)|\le \Lambda_n(x) \epsilon, \qquad \epsilon= \max_{0\le k\le n}|f(x_k)-F(x_k)| .
\]
Moreover, using the polynomial reproducing property (\ref{eq-inva}), it is easy to prove that the Lebesgue constants are also involved in the error estimate
\[
|\tilde r(f,x)-f(x)|\le [1 + \Lambda_n] E_d(f),\qquad x\in [a,b]
\]
where $ E_d(f)$ denotes the error of best approximation of $f$ in $\PP_d$ w.r.t. the uniform norm $\|f\|_\infty=\sup_{x\in [a,b]}|f(x)|$, namely
\begin{equation}\label{E-d}
E_d(f)=\inf_{P\in\PP_d} \|f-P\|_\infty.
\end{equation}
For classical FH interpolation ($\gamma=1$) the Lebesgue constants have been estimated in \cite{bos2011lebesgue,bos2012lebesgue} for equidistant nodes and in \cite{hormann2012barycentric} for quasi--equidistant nodes. In both cases, they result to grow as $\log n$ with $n$ and as $2^d$ with $d$.
Here we show that, introducing the additional parameter $\gamma>1$, for the generalized FH interpolation we succeed in  getting Lebesgue constants uniformly bounded w.r.t. $n$.

More precisely, setting
\begin{equation}\label{h}
h=\max_{0\le k<n}(x_{k+1}-x_k),\qquad h^*=\min_{0\le k<n}(x_{k+1}-x_k)
\end{equation}
we have the following
\begin{theorem}\label{th-LC}
For all $n,d,\gamma\in\NN$, with $d\in [1,n]$ and $\gamma>1$, we have
\begin{equation}\label{eq-LC}
\Lambda_n(x)\le \C\ d 2^d \left(\frac h{h^*}\right)^{\gamma+d}, \qquad \forall x\in [a,b],
\end{equation}
where $\C>0$ is a constant independent of $n,h,x,d$ and $\gamma$.
\end{theorem}
\begin{remark}\label{rem-LC}
 Similarly to the classic FH interpolation, we note that the dependence on $d$ of the Lebesgue constants is exponential for all $\gamma>1$ too. Indeed, we conjecture the linear factor $d$ in (\ref{eq-LC}) can be removed, but we were not able to prove it.
\end{remark}
An immediate consequence of Thm. \ref{th-LC} is the following
\begin{corollary}\label{cor-lc}
Let $n,d,\gamma\in\NN$, with $d\in [1,n]$ and $\gamma>1$. In the case of equidistant nodes (i.e., if  $h=h^*$) and, more generally,  in the case of quasi--equidistant nodes (i.e., if  $h/h^*\le\rho$ holds with $\rho\in\RR$ independent of $n$), we get
\[
\sup_n\Lambda_n <\infty .
\]
\end{corollary}
\subsection{Proof of Theorem \ref{th-LC}}\label{sec:LC2}
In order to prove Thm. \ref{th-LC} let us first state three preliminary lemmas
\begin{lemma}\label{lem-LC}
Let  $n,d,\gamma\in\NN$ with $1\le d\le n$ and $\gamma>1$. If $x_\ell<x<x_{\ell +1}$, with $0\le\ell<n$, then for all $i\in\{0,\ldots, n-d\}$ and any $k\in\{i,\ldots, i+d\}$, we have
\begin{equation}\label{eq-lem-LC}
\frac 1{|D(x)|}\prod_{s=i,\ s\ne k}^{i+d}\frac 1{|x-x_s|^{\gamma-1}}\le \left\{\begin{array}{lll}
\displaystyle |x-x_{\ell+1}|^\gamma \prod_{s=\ell -d+1}^\ell{|x-x_s|} & \mbox{if}& i\le \ell-d\\ [.2in]
\displaystyle |x-x_k|^\gamma\prod_{s=i, s\ne k}^{i+d}{|x-x_s|} & \mbox{if}& i\in I_\ell\\ [.2in]
\displaystyle  |x-x_{\ell}|^\gamma\prod_{s=\ell+1}^{\ell +d}{|x-x_s|} & \mbox{if}& i\ge \ell+1
\end{array} \right.
\end{equation}
where $D(x)$ and $I_\ell$ are defined in (\ref{D}) and (\ref{Il}), respectively.
\end{lemma}
{\it Proof of Lemma \ref{lem-LC}}\newline
Firstly note that the following inequality
\begin{equation}\label{eq-prod-LC}
\prod_{s=i,\ s\ne k}^{i+d}\frac 1{|x-x_s|^{\gamma-1}}\le \left\{\begin{array}{lll}
\displaystyle \prod_{s=\ell -d+1}^\ell\frac 1{|x-x_s|^{\gamma-1}} & \mbox{if}& i\le \ell-d\\ [.2in]
\displaystyle \prod_{s=\ell+1}^{\ell +d}\frac 1{|x-x_s|^{\gamma-1}} & \mbox{if}& i\ge \ell+1
\end{array} \right.
\end{equation}
can be proved by taking into account that
\[
|x-x_{i+j}|\ge \left\{\begin{array}{lll}
\displaystyle |x-x_{\ell-d+1+j}| & \mbox{if}& 0\le j< k-i\\ [.1in]
\displaystyle |x-x_{\ell-d+j}| & \mbox{if}& k-i<j\le d
\end{array} \right.\qquad \forall i\le \ell-d
\]
\[
|x-x_{i+j}|\ge \left\{\begin{array}{lll}
\displaystyle |x-x_{\ell+1+j}| & \mbox{if}& 0\le j< k-i\\ [.1in]
\displaystyle |x-x_{\ell+j}| & \mbox{if}& k-i<j\le d
\end{array} \right.\qquad \forall i\ge \ell+1
\]
In addition, by (\ref{rem-lambda}), we deduce
\begin{equation}\label{eq-D-LC}
|D(x)| \ge \left\{\begin{array}{lll}
\displaystyle |\tlambda_{\ell-d+1}(x)|=\prod_{s=\ell -d+1}^{\ell+1}\frac 1{|x-x_s|^{\gamma}} & \mbox{if}& i\le \ell-d\\ [.2in]
\displaystyle |\tlambda_{i}(x)|=\prod_{s=i}^{i+d}\frac 1{|x-x_s|^{\gamma}} & \mbox{if}& \ell -d+1\le i\le \ell\\ [.2in]
\displaystyle \displaystyle |\tlambda_{\ell}(x)|=\prod_{s=\ell}^{\ell +d}\frac 1{|x-x_s|^{\gamma}} & \mbox{if}& i\ge \ell+1
\end{array} \right.
\end{equation}
Hence,  by collecting (\ref{eq-prod-LC}) and (\ref{eq-D-LC}), we obtain (\ref{eq-lem-LC}).
\Proofend

\begin{lemma}\label{lem-LC1}
Let $n,d,\gamma\in\NN$ with $1\le d\le n$ and $\gamma>1$. If $x_\ell<x<x_{\ell +1}$, with $0\le\ell<n$, then for all $i\in\{0,\ldots, n-d\}$ and any $k\in\{i,\ldots, i+d\}$, we have
\begin{equation}\label{eq-lem-LC1}
\frac 1{|D(x)|}\prod_{s=i,\ s\ne k}^{i+d}\frac 1{|x-x_s|^{\gamma-1}}\le \left\{\begin{array}{lll}
 h^{d+\gamma}d! & \mbox{if}& i \in\bar{I}_\ell\\ [.2in]
\displaystyle h^d|x-x_k|^\gamma (\ell+1-i)!(d+i-\ell)!& \mbox{if}& i\in I_\ell
\end{array} \right.
\end{equation}
where $\overline{I}_\ell=\{i\in\{0,\ldots, n-d\} : \ i\le \ell -d \ \vee\ i\ge \ell+1\}$ is the  complementary set of $I_\ell$ defined in (\ref{Il}).
\end{lemma}
{\it Proof of Lemma \ref{lem-LC1}}\newline
Taking into account that
\begin{equation}\label{equi}
|x_k-x_s|\le h|k-s|,\qquad \forall k,s\in\{0,\ldots, n\}
\end{equation}
we get
\begin{eqnarray}
\label{prod1}
\prod_{s=l-d+1}^{\ell}|x-x_s|&\le& \prod_{s=l-d+1}^{\ell}|x_{\ell+1}-x_s|\le h^d\prod_{s=l-d+1}^{\ell}|l+1-s|=h^d d!\\
\label{prod2}
\prod_{s=l+1}^{\ell+d}|x-x_s|&\le& \prod_{s=l+1}^{\ell+d}|x_{\ell}-x_s|\le h^d\prod_{s=l+1}^{\ell+d}(\ell-s)=h^d d!
\end{eqnarray}
Also, by (\ref{equi}), $\forall i\in I_\ell$ we have
\begin{eqnarray}
\nonumber
\prod_{s=i, \ s\ne k}^{i+d}|x-x_s|&\le&
\prod_{s=i, \ s\ne k}^{\ell}|x_{\ell+1}-x_s|\prod_{s=\ell+1, \ s\ne k}^{i+d}|x_{\ell}-x_s|\\
\label{prod3}
&\le& h^d (\ell+1-i)! (d+i-\ell)!
\end{eqnarray}
Hence, the statement follows by applying (\ref{prod1})--(\ref{prod3}) to the result of Lemma \ref{lem-LC}. \Proofend
\begin{lemma}\label{lem-LC2}
For all $n\in\NN$ and $i\in J_k=\{i\in\{0,\ldots, (n-d)\} \ : \ k-d\le i\le k\}$, with $k=0,\ldots, n$, we have
\begin{equation}\label{eq-lem-LC2}
\prod_{s=i, \ s\ne k}^{i+d}\frac 1 {|x_k-x_s|}\le \left(\frac 1{h^*}\right)^d\frac 1{(k-i)! (d+i-k)!} 
\end{equation}
\end{lemma}
{\it Proof of Lemma \ref{lem-LC2}}\newline
Since
\begin{equation}\label{equi1}
|x_k-x_s| \ge h^*|k-s|,\qquad \forall k, s\in\{0,\ldots,n\}
\end{equation}
we get
\begin{equation}\label{eq-Aik}
\prod_{s=i, \ s\ne k}^{i+d}\frac 1 {|x_k-x_s|}\le \left(\frac 1{h^*}\right)^d
\prod_{s=i, \ s\ne k}^{i+d}\frac 1 {|k-s|}= \left(\frac 1{h^*}\right)^d\frac 1{(k-i)! (d+i-k)!}
\end{equation}
\Proofend

{\bf Proof of Theorem \ref{th-LC}}

First of all, note that if $x\in \{x_0,\ldots, x_n\}$ then the statement is trivial since,  by (\ref{bk}), we have
\[
\Lambda_n(x_k)=1, \qquad k=0,1,\ldots,n.
\]
Hence, let us assume $x_\ell<x<x_{\ell +1}$ with $\ell=0,\ldots, (n-1)$, and prove
(\ref{eq-LC}).

By (\ref{bk}) and (\ref{wk}) we get
\begin{eqnarray*}
\Lambda_n(x)&=&\sum_{k=0}^n|b_k(x)|=\sum_{k=0}^n
\frac{|w_k(x)|}{|x-x_k|^\gamma |D(x)|} \\
&\le& \sum_{k=0}^n\sum_{i\in J_k}\frac 1{|x-x_k|^\gamma|D(x)|}
\prod_{s=i, \ s\ne k}^{i+d}\frac 1 {|x_k-x_s| |x-x_s|^{\gamma-1}}
\end{eqnarray*}
i.e., set for brevity
\begin{equation}\label{ABik}
A_{i,k}=\prod_{s=i, \ s\ne k}^{i+d}\frac 1 {|x_k-x_s| },\qquad
B_{i,k}(x)=\prod_{s=i, \ s\ne k}^{i+d}\frac 1 {|x-x_s|^{\gamma-1}}
\end{equation}
we have
\begin{equation}\label{sum}
\Lambda_n(x)\le \sum_{k=0}^n\sum_{i\in J_k}\frac {A_{i,k}\ B_{i,k}(x)}{|x-x_k|^\gamma|D(x)|} 
\end{equation}
Now we note that
\[
J_k=(J_k \cap \overline{I}_\ell)\cup (J_k\cap I_\ell), \qquad k=0.\ldots, n,
\]
where $I_\ell$ is given by (\ref{Il}) and $\overline{I}_\ell$ is the complementary set.

Hence, we decompose the summation in (\ref{sum})
\begin{eqnarray*}
\Lambda_n(x)
&\le&  \sum_{k=0}^n\left\{\sum_{i\in J_k\cap \overline{I}_\ell}
\frac{A_{i,k}\ B_{i,k}(x)}{|x-x_k|^\gamma|D(x)|}+ \sum_{i\in J_k\cap I_\ell}
\frac{A_{i,k}\ B_{i,k}(x)}{|x-x_k|^\gamma|D(x)|}\right\}
\\
&=:& \sum_{k=0}^n \left\{\sigma_k(x) + \mu_k(x)\right\}=:S_1(x)+S_2(x)
\end{eqnarray*}
where, as usual, we mean that $\sum_{i\in I}a_i=0$ whenever $I$ is the empty set.

Let us first estimate
\[
S_1(x)= \sum_{k=0}^n\sigma_k(x)=\sum_{k=0}^n\left(\sum_{i\in J_k\cap \overline{I}_\ell}\frac {A_{i,k}\ B_{i,k}(x)}{|x-x_k|^\gamma|D(x)|}\right)
\]
by distinguishing the following cases:
\begin{itemize}
\item \textsc{Case $k=\ell$}. Let us estimate the term $\sigma_\ell(x)$

Note that when $k=\ell$ we have $J_k\cap \overline{I}_\ell=\{\ell-d\}$ if $\ell\ge d$, otherwise it is $J_k\cap \overline{I}_\ell=\emptyset$. Hence, suppose $\ell\ge d$ (otherwise $\sigma_\ell(x)=0$  holds), we have
\[
\sigma_\ell(x):=\sum_{i\in J_\ell\cap \overline{I}_\ell} \frac{A_{i,\ell}\ B_{i,\ell}(x)}{|x-x_\ell|^\gamma|D(x)|}=
\frac{A_{\ell-d,\ell}\ B_{\ell-d,\ell}(x)}{|x-x_\ell|^\gamma|D(x)|}
\]
In this case, we use (\ref{rem-lambda}) as follows
\[
|D(x)|\ge |\tlambda_{\ell-d+1}|=\prod_{s=\ell-d+1}^{\ell +1}\frac 1{|x-x_s|^\gamma}
\]
and by Lemma \ref{lem-LC2}, we have
\[
A_{\ell-d,\ell}\le \left(\frac 1{h^*}\right)^d\frac 1{d!}
\]
Consequently, we get
    \begin{eqnarray*}
\sigma_\ell(x)&=&\frac{A_{\ell-d,\ell}\ B_{\ell-d,\ell}(x)}{|x-x_\ell|^\gamma|D(x)|}\\
&\le& \frac 1{d!(h^*)^d}\frac 1{|x-x_\ell|^\gamma |\tlambda_{\ell-d+1}(x)|}
\prod_{s=\ell-d}^{\ell-1}\frac 1 {|x-x_s|^{\gamma-1}}\\
&=&\frac 1{d!(h^*)^d}\frac{|x-x_{\ell+1}|^\gamma}{|x-x_{\ell-d}|^{\gamma}}
\prod_{s=\ell-d}^{\ell-1}|x-x_s|\\
&\le& \frac 1{d!(h^*)^d}
\frac{|x_\ell-x_{\ell+1}|^\gamma}{|x_{\ell-d+1}-x_{\ell-d}|^{\gamma}}
\prod_{s=\ell-d}^{\ell-1}|x_{\ell +1}-x_s|
    \end{eqnarray*}
    and by  (\ref{equi}), (\ref{equi1}), we conclude
\begin{equation}\label{kl}
\sigma_\ell(x)\le \left(\frac h{h^*}\right)^{d+\gamma}\frac 1{d!}\prod_{s=\ell-d}^{\ell-1}|\ell+1-s|
= \left(\frac h{h^*}\right)^{d+\gamma}(d+1) .
\end{equation}
\item \textsc{Case $k=\ell+1$}. Let us estimate the term $\sigma_{\ell+1}(x)$

When $k=\ell+1$ we have $J_k\cap \overline{I}_\ell=\{\ell+1\}$ and $\sigma_{\ell+1}(x)$ can be estimated by means of (\ref{rem-lambda}), (\ref{eq-lem-LC2}), (\ref{equi}), and (\ref{equi1}) similarly to the previous case, getting
\begin{eqnarray}
\nonumber
\sigma_{\ell+1}(x)&:=&\sum_{i\in J_{\ell+1}\cap \overline{I}_\ell} \frac{A_{i,\ell+1}\ B_{i,\ell+1}(x)}{|x-x_{\ell+1}|^\gamma|D(x)|}=
\frac{A_{\ell+1,\ell+1}\ B_{\ell+1,\ell+1}(x)}{|x-x_{\ell+1}|^\gamma|D(x)|}\\
\nonumber
&\le& \frac 1{d!(h^*)^d}\frac 1{|x-x_{\ell+1}|^\gamma |\tlambda_{\ell}(x)|}
\prod_{s=\ell+2}^{\ell+1+d}\frac 1 {|x-x_s|^{\gamma-1}}\\
\nonumber
&=&\frac 1{d!(h^*)^d}\frac{|x-x_{\ell}|^\gamma}{|x-x_{\ell+1+d}|^{\gamma}}
\prod_{s=\ell+2}^{\ell+1+d}|x-x_s|\\
\nonumber
&\le&\frac 1{d!(h^*)^d}
\frac{|x_\ell-x_{\ell+1}|^\gamma}{|x_{\ell+d}-x_{\ell+1+d}|^{\gamma}}
\prod_{s=\ell+2}^{\ell+1+d}|x_\ell-x_s|\\
\label{kl1}
&\le& \left(\frac h{h^*}\right)^{d+\gamma}\frac 1{d!}\prod_{s=\ell+2}^{\ell+1+d}|\ell-s|
= \left(\frac h{h^*}\right)^{d+\gamma} (d+1).
    \end{eqnarray}
\item \textsc{Case $k\not\in\{\ell,\ell+1\}$}. Let us estimate the summation of the remaining terms $\sigma_k(x)$. 

By applying Lemma \ref{lem-LC1} and Lemma \ref{lem-LC2} we have
\begin{eqnarray*}
\sum_{\scriptsize\begin{array}{c}k=0\\ k\not\in\{\ell,\ell+1\}\end{array}}^n \sigma_k(x)&=&
\sum_{\scriptsize\begin{array}{c}k=0\\ k\not\in\{\ell,\ell+1\}\end{array}}^n
\sum_{i\in J_k\cap \overline{I}_\ell}
\left(\frac {A_{i,k}\ B_{i,k}(x)}{|x-x_k|^\gamma|D(x)|}\right)\\
&\le& \sum_{\scriptsize\begin{array}{c}k=0\\ k\not\in\{\ell,\ell+1\}\end{array}}^n \sum_{i\in J_k\cap \overline{I}_\ell}
\frac{h^\gamma}{|x-x_k|^\gamma}\left(\frac h{h^*}\right)^d
\left(\begin{array}{c}d\\ k-i\end{array}\right)
\end{eqnarray*}
and taking into account that
\[
\sum_{i\in J_k\cap \overline{I}_\ell}\left(\begin{array}{c}d\\ k-i\end{array}\right)\le \sum_{i\in J_k} \left(\begin{array}{c}d\\ k-i\end{array}\right)\le \sum_{j=0}^d \left(\begin{array}{c}d\\ j\end{array}\right)=2^d
\]
we continue the estimate as follows
\begin{eqnarray}
\nonumber
\sum_{\scriptsize\begin{array}{c}k=0\\ k\not\in\{\ell,\ell+1\}\end{array}}^n \sigma_k(x)
&\le& h^\gamma \left(\frac h{h^*}\right)^d  \sum_{\scriptsize\begin{array}{c}k=0\\ k\not\in\{\ell,\ell+1\}\end{array}}^n \frac 1{|x-x_k|^\gamma}\sum_{i\in J_k\cap \overline{I}_\ell}\left(\begin{array}{c}d\\ k-i\end{array}\right)\\
\nonumber
&\le&2^d\ h^\gamma \left(\frac h{h^*}\right)^d \left(
\sum_{k=0}^{\ell-1}\frac 1{|x_\ell-x_k|^\gamma}
+\sum_{k=\ell+2}^n \frac 1{|x_{\ell+1}-x_k|^\gamma}\right) \\
\nonumber
&\le& 2^d \left(\frac h{h^*}\right)^{\gamma+d} \left(
\sum_{k=0}^{\ell-1}\frac 1{|\ell-k|^\gamma}
+\sum_{k=\ell+2}^n \frac 1{|{\ell+1}-k|^\gamma}\right)\\
\label{kl2}
&\le & 2^{d+1}\left(\frac h{h^*}\right)^{\gamma+d} \sum_{j=1}^n\frac 1{j^\gamma}\le \C 2^d \left(\frac h{h^*}\right)^{\gamma+d} .
\end{eqnarray}
having used $\sum_{j=1}^n \frac 1{j^\gamma}\le \sum_{j=1}^\infty \frac 1{j^2} <\infty$ in the last inequality.
\end{itemize}
Summing up, by (\ref{kl})--(\ref{kl2}) we have proved that $S_1(x)$ satisfies the bound in (\ref{eq-LC}).

Now let us prove that the same holds for
\[
S_2(x)= \sum_{k=0}^n\mu_k(x)=\sum_{k=0}^n\left( \sum_{i\in J_k\cap I_\ell} \frac {A_{i,k}\ B_{i,k}(x)}{|x-x_k|^\gamma|D(x)|}\right) .
\]
We note that if $k\le \ell-d$ or $k>\ell+d$ then we certainly have $J_k\cap I_\ell=\emptyset$. Thus, $S_2(x)$ is, indeed, given by the following sum
\[
S_2(x)= \sum_{\scriptsize\begin{array}{c}k=\ell +1-d\\ [-.08cm]
0\le k\le n\end{array}}^{\ell +d}\left(\sum_{i\in J_k\cap I_\ell}\frac {A_{i,k}\ B_{i,k}(x)}{|x-x_k|^\gamma|D(x)|}\right)
\]
Hence,  by using Lemma \ref{lem-LC1} and Lemma \ref{lem-LC2},  we get
\begin{eqnarray*}
S_2(x)
&\le& \left(\frac h{h^*}\right)^d\sum_{\scriptsize\begin{array}{c}k=\ell +1-d\\ [-.08cm]
0\le k\le n\end{array}}^{\ell +d}\sum_{i\in J_k\cap I_\ell}
\frac{(\ell+1-i)!(i+d-\ell)!}{(k-i)!(i+d-k)!} .
\end{eqnarray*}
Now we distinguish the following cases:
\begin{itemize}
\item If $k< \ell$ then we have $(\ell+1-i)>(k-i)$ and $(i+d-\ell)<(i+d-k)$. Hence we can write
    \[
     \frac{(\ell+1-i)!(i+d-\ell)!}{(k-i)!(i+d-k)!}= \frac{(\ell+1-i)(\ell+1-i-1)\cdots (k-i+1)}{(i+d-k)(i+d-k-1)\cdots (i+d-\ell +1))}
    \]
    and taking into account that, for $i\in I_\ell$, the right hand side term takes its maximum value when $i=\ell+1-d$, we get
    \begin{equation}\label{case1}
    \frac{(\ell+1-i)!(i+d-\ell)!}{(k-i)!(i+d-k)!}\le\frac{d (d-1)\cdots (d-(\ell -k))}{(\ell+1-k)(\ell-k)\cdots 2}=\left(\begin{array}{c}d\\ \ell+1-k \end{array}\right) .
    \end{equation}
    \item If $k=\ell$ then we have
    \begin{equation}\label{case2}
    \frac{(\ell+1-i)!(i+d-\ell)!}{(k-i)!(i+d-k)!}=(\ell+1-i)\le d .
    \end{equation}
    \item If $k=\ell+1$ then we have
    \begin{equation}\label{case3}
    \frac{(\ell+1-i)!(i+d-\ell)!}{(k-i)!(i+d-k)!}=(i+d-\ell)\le d .
    \end{equation}
 \item If $k> \ell +1$ then we have $(\ell+1-i)<(k-i)$ and $(i+d-\ell)>(i+d-k)$. Hence we can write
     \[
     \frac{(\ell+1-i)!(i+d-\ell)!}{(k-i)!(i+d-k)!}= \frac{(i+d-\ell)(i+d-\ell-1)\cdots (i+d-k +1))}{(k-i)(k-i-1)\cdots (\ell +2-i)}
    \]
    and since, for $i\in I_\ell$, the maximum value of the right hand side term is achieved when $i=\ell$, we get
    \begin{equation}\label{case4}
    \frac{(\ell+1-i)!(i+d-\ell)!}{(k-i)!(i+d-k)!}\le\frac{d (d-1)\cdots (d-(k-\ell-1))}{(k-\ell)(k-\ell-1)\cdots 2}=\left(\begin{array}{c}d\\ k-\ell \end{array}\right) .
    \end{equation}
\end{itemize}
Thus, by virtue of (\ref{case1})--(\ref{case4}), we have proved
\[
\frac{(\ell+1-i)!(i+d-\ell)!}{(k-i)!(i+d-k)!}\le \left\{\begin{array}{ll}
\left(\begin{array}{c}d\\ \ell+1-k \end{array}\right) & \mbox{if $\ \ell+1-d\le k\le \ell$}\\ [.15in]
\left(\begin{array}{c}d\\ k-\ell \end{array}\right) & \mbox{if $\ \ell+1\le k\le \ell+d$}
\end{array}\right.
\]

Consequently, taking into account that the set $J_k\cap I_\ell$ has at most $d$ elements, the estimate of $S_2(x)$ continues as follows
\begin{eqnarray}
\nonumber
S_2(x)&\le& \left(\frac h{h^*}\right)^d
\sum_{\scriptsize\begin{array}{c}k=\ell +1-d\\ [-.08cm]
0\le k\le n\end{array}}^{\ell +d}\
\sum_{i\in J_k\cap I_\ell}
\frac{(\ell+1-i)!(i+d-\ell)!}{(k-i)!(i+d-k)!}\\
\nonumber
&\le& \left(\frac h{h^*}\right)^d
\left[\sum_{\scriptsize\begin{array}{c}k=\ell +1-d\\ [-.08cm]
0\le k\le n\end{array}}^{\ell}\
\sum_{i\in J_k\cap I_\ell}\left(\begin{array}{c}d\\ \ell+1-k \end{array}\right)+
\sum_{k=\ell +1}^{\ell +d}\sum_{i\in J_k\cap I_\ell}\left(\begin{array}{c}d\\ k-\ell \end{array}\right)
\right]\\
\label{eq-S1}
&\le& \left(\frac h{h^*}\right)^d  d\ 2 \sum_{j=1}^d\left(\begin{array}{c}d\\ j \end{array}\right)=2 \left(\frac h{h^*}\right)^d  d \ 2^d .
\end{eqnarray}
\subsection{A related result}\label{sec:LC3}
Going along the same lines as the proof of Thm. \ref{th-LC}, we get the following result which will be useful in the next section.
\begin{theorem}\label{th-LCrelated}
For all $n,\gamma,d\in\NN$, with $1\le d\le n$ and $\gamma>1$,  let $\{b_k(x)\}_k$ be the interpolating basis defined in (34). If the distribution of nodes is such that $h\sim h^*\sim n^{-1}$ holds ($a\sim b$ meaning that the ratio $a/b$ is between two absolute, positive, constants) then for all $\alpha>0$ we have
\begin{equation}\label{eq-LCrelated}
\Sigma(x):=\sum_{k=0}^n |x-x_k|^\alpha |b_k(x)|\le \C \left\{\begin{array}{lll}
\displaystyle\frac 1{n^\alpha} & \mbox{if} & \gamma >1+\alpha\\ [.15in]
\displaystyle\frac{\log n} {n^\alpha} & \mbox{if} & \gamma=1+\alpha\\ [.15in]
\displaystyle\frac 1{n^{\gamma-1}} & \mbox{if} & 1<\gamma<1+\alpha
\end{array}\right. \qquad \forall x\in [a,b],
\end{equation}
where $\C>0$ is a constant independent of $n, x$, bounded with respect to $\gamma>1$ but exponentially growing with $d$.
\end{theorem}
{\it Proof of Theorem \ref{th-LCrelated}}.
Recalling that $b_k(x_j)=\delta_{k,j}$ (cf. (\ref{bk})),  if $x\in\{x_0,\ldots, x_n\}$ then the statement is trivial since $\Sigma(x)=0$. Hence, let us fix $x_\ell<x<x_{\ell+1}$ with $\ell=0,\ldots,n-1$.
Following the same lines as in the proof of Thm. \ref{th-LC} and using the notations therein introduced, we note that
\begin{eqnarray*}
\Sigma(x)&:=&\sum_{k=0}^n |x-x_k|^\alpha |b_k(x)|\\
&\le&
\sum_{k=0}^n |x-x_k|^\alpha\sum_{i\in J_k}\frac{A_{i,k}\ B_{i,k}(x)}{|x-x_k|^\gamma |D(x)|}\\
&=&\sum_{k=0}^n |x-x_k|^\alpha\left\{\sum_{i\in J_k\cap \overline{I}_\ell}
\frac{A_{i,k}\ B_{i,k}(x)}{|x-x_k|^\gamma|D(x)|}+ \sum_{i\in J_k\cap I_\ell}
\frac{A_{i,k}\ B_{i,k}(x)}{|x-x_k|^\gamma|D(x)|}\right\}\\
&=& \sum_{k=0}^n |x-x_k|^\alpha\sigma_k(x) +\sum_{k=0}^n |x-x_k|^\alpha\mu_k(x)=:\Sigma_1(x)+\Sigma_2(x) .
\end{eqnarray*}
Hence, we estimate $\Sigma_1(x)$ and $\Sigma_2(x)$ using the results obtained in the proof of Thm. \ref{th-LC} for $S_1(x)$ and $S_2(x)$, respectively.

For simplicity, in the sequel we denote by $\C$ all positive constants as in the statement, even if they have different values.

As regards $\Sigma_1(x)$, similarly to $S_1(x)$ in the proof of Thm. \ref{th-LC}, we deduce
\begin{eqnarray*}
\Sigma_1(x)&=& |x-x_\ell|^\alpha \sigma_\ell(x)+ |x-x_{\ell+1}|^\alpha\sigma_{\ell+1}(x)+
\sum_{\scriptsize\begin{array}{c}k=0\\ k\not\in\{\ell,\ell+1\}\end{array}}^n
|x-x_k|^\alpha\sigma_k(x)\\
&\le& 2(d+1)\left(\frac h{h^*}\right)^d h^\alpha +
\left(\frac h{h^*}\right)^d h^\gamma\sum_{\scriptsize\begin{array}{c}k=0\\ k\not\in\{\ell,\ell+1\}\end{array}}^n
\sum_{i\in J_k\cap \overline{I}_\ell}
\frac{|x-x_k|^\alpha}{|x-x_k|^\gamma}
\left(\begin{array}{c}d\\ k-i\end{array}\right)\\
&\le& 2(d+1)\left(\frac h{h^*}\right)^d h^\alpha +
2^d \left(\frac h{h^*}\right)^d h^\gamma\left(
\sum_{k=0}^{\ell-1}\frac {|x_{\ell +1}-x_k|^\alpha}{|x_\ell-x_k|^\gamma}
+\sum_{k=\ell+2}^n \frac {|x_\ell-x_k|^\alpha}{|x_{\ell+1}-x_k|^\gamma}\right)
\\
&\le& 2(d+1)\left(\frac h{h^*}\right)^d h^\alpha +
2^d \left(\frac h{h^*}\right)^{d+\gamma} h^\alpha\left(
\sum_{k=0}^{\ell-1}\frac {|\ell +1-k|^\alpha}{|\ell -k|^\gamma}
+\sum_{k=\ell+2}^n \frac {|\ell -k|^\alpha}{|{\ell+1}-k|^\gamma}\right)\\
&\le& \C  \left(\frac h{h^*}\right)^{d+\gamma} h^\alpha\left(
\sum_{j=1}^{n}\frac 1{j^{\gamma-\alpha}}\right)\le
\C\frac 1{n^\alpha}\left(\sum_{j=1}^{n}\frac 1{j^{\gamma-\alpha}}\right)
\end{eqnarray*}
having used the hypothesis $h\sim h^*\sim n^{-1}$ to get the last inequality.

Now we observe that for all $p>0$ it is
\[
\sum_{j=1}^n\frac 1{j^p}\le 1 +\sum_{j=1}^n\frac 1{(j+1)^p}\le 1+\sum_{j=1}^n\int_j^{j+1}\frac{dx}{x^p}=1+\int_1^{(n+1)}\frac{dx}{x^p}=
\left\{\begin{array}{lll}
 1+\log (n+1) & \mbox{if} & p=1\\ [.1in]
\frac{(n+1)^{1-p}-p} {1-p} & \mbox{if} & p\ne 1
\end{array}\right.
\]
and for all $p<0$ it is
\[
\sum_{j=1}^n\frac 1{j^p}\le \sum_{j=1}^n\int_j^{j+1}x^{-p}dx=\int_1^{n+1}x^{-p}dx=\frac{(n+1)^{1-p}-1}{1-p}.
\]
Hence, the following holds
\begin{equation}\label{armonic-sum}
\sum_{j=1}^{\infty}\frac 1{j^p}\le  \left\{\begin{array}{lll}
\C \log n & \mbox{if} & p=1\\ [.1in]
\C n^{1-p} & \mbox{if} & p< 1\\ [.1in]
\frac p{p-1} & \mbox{if} & p> 1
\end{array}\right., \qquad \forall p\in\RR
\end{equation}
which implies
\[
\Sigma_1(x) \le \frac \C {n^\alpha}\left(\sum_{j=1}^{n}\frac 1{j^{\gamma-\alpha}}\right) \le  \frac \C {n^\alpha}
 \left\{\begin{array}{lll}
 \log n & \mbox{if} & \gamma-\alpha=1\\ [.1in]
 n^{1-\gamma+\alpha} & \mbox{if} &\gamma-\alpha< 1\\ [.1in]
1 & \mbox{if} & \gamma-\alpha> 1
\end{array}\right.
\]
Finally, as regards $\Sigma_2(x)$, we deduce the following from the results achieved for $S_2(x)$ in the proof of Thm. \ref{th-LC}
\begin{eqnarray*}
\Sigma_2(x)&=&\sum_{k=0}^n|x-x_k|^\alpha \left( \sum_{i\in J_k\cap I_\ell} \frac {A_{i,k}\ B_{i,k}(x)}{|x-x_k|^\gamma|D(x)|}\right)\\
&=& \sum_{\scriptsize\begin{array}{c}k=\ell +1-d\\ [-.08cm]
0\le k\le n\end{array}}^{\ell +d}|x-x_k|^\alpha
\left(\sum_{i\in J_k\cap I_\ell}
\frac {A_{i,k}\ B_{i,k}(x)}{|x-x_k|^\gamma|D(x)|}\right)
\\
&\le& \left(\frac h{h^*}\right)^d \sum_{\scriptsize\begin{array}{c}k=\ell +1-d\\ [-.08cm]
0\le k\le n\end{array}}^{\ell +d}\ |x-x_k|^\alpha \sum_{i\in J_k\cap I_\ell}
\frac{(\ell+1-i)!(i+d-\ell)!}{(k-i)!(i+d-k)!}\\
&\le& d \left(\frac h{h^*}\right)^d
\left[\sum_{\scriptsize\begin{array}{c}k=\ell +1-d\\ [-.08cm]
k\ge 0\end{array}}^{\ell}\ |x_{\ell+1}-x_k|^\alpha\left(\begin{array}{c}d\\ \ell+1-k \end{array}\right)
+
\sum_{\scriptsize\begin{array}{c}k=\ell +1\\ [-.08cm]
 k\le n\end{array}}^{\ell +d}|x_\ell-x_k|^\alpha\left(\begin{array}{c}d\\ k-\ell \end{array}\right)
\right]
\\
&\le& d\ h^\alpha\left(\frac h{h^*}\right)^d
\left[\sum_{\scriptsize\begin{array}{c}k=\ell +1-d\\ [-.08cm]
0\le k\le n\end{array}}^{\ell}\ |{\ell+1}-k|^\alpha
\left(\begin{array}{c}d\\ \ell+1-k \end{array}\right)+
\sum_{\scriptsize\begin{array}{c}k=\ell +1\\ [-.08cm]
0\le k\le n\end{array}}^{\ell +d}|\ell-k|^\alpha
\left(\begin{array}{c}d\\ k-\ell \end{array}\right)
\right]
\\
&\le& 2d\ h^\alpha\left(\frac h{h^*}\right)^d\left[\sum_{j=1}^d
j^\alpha\left(\begin{array}{c}d\\ j \end{array}\right)
\right]\le \C h^\alpha\left(\frac h{h^*}\right)^d \le \frac \C{n^\alpha}
\end{eqnarray*}
having used $h\sim h^*\sim n^{-1}$ to get the last inequality.
\Proofend
\section{Error estimates}\label{sec:er}
First of all let us state the following fundamental result
\begin{theorem}\label{th-conv}
 Let the parameters  $d,\gamma\in\NN$ be arbitrarily fixed with $\gamma\ge 2$. Moreover, let the distribution of nodes satisfy $h\sim h^*\sim n^{-1}$. as $n\to\infty$. Then for any continuous function $f\in C([a,b])$ we have
\begin{equation}\label{eq-conv}
\lim_{n\to \infty}\|\tilde r(f)-f\|_\infty=0
\end{equation}
\end{theorem}
{\it Proof of Theorem \ref{th-conv}}.
Let us arbitrarily fix $d,\gamma\in\NN$ with $\gamma\ge 2$, and consider arbitrarily large integers $n\ge d$. First, we prove that (\ref{eq-conv}) holds if $f=P$ is any polynomial. Due to (\ref{eq-inva}), this is trivial if $\deg P\le d$. Hence let $P$ be a polynomial of degree $s>d\ge 1$

Taking into account that $\tilde r(f, x)$ certainly preserves the constant function $f(x)=1$ (since $d\ge 1$), by (\ref{r-bk}) we have
\[
1=\sum_{k=0}^n b_k(x), \qquad \forall x\in [a,b],
\]
and consequently, for all integers $n\ge d$, we get
\begin{equation}\label{f-rf}
f(x)-\tilde r(f,x)=\sum_{k=0}^n \left[f(x)-f(x_k)\right] b_k(x), \qquad \forall x\in [a,b], \quad \forall f\in C([a,b]).
\end{equation}
On the other hand, by the Mean Value Theorem, we note that
\begin{equation}\label{lip-P}
|P(x)-P(y)|\le M |x-y|, \qquad \forall x,y\in [a,b], \qquad M=\|P^\prime\|_\infty>0.
\end{equation}
In conclusion, by collecting (\ref{f-rf})  and (\ref{lip-P}),  $\forall x\in [a,b]$, we get
\[
|P(x)-\tilde r(P,x)|
 \le \sum_{k=0}^n \left|P(x)-P(x_k)\right| |b_k(x)|
\le  M  \sum_{k=0}^n \left|x-x_k\right||b_k(x)|,
\]
and applying Thm. \ref{th-LCrelated} we obtain
\[
\|P-\tilde r(P)\|_\infty \le\C\left \{\begin{array}{lll}
\displaystyle\frac 1{n} & \mbox{if} & \gamma >2\\ [.15in]
\displaystyle\frac{\log n} {n} & \mbox{if} & \gamma=2
\end{array}\right.
\]
where $\C>0$ is independent of $n$.\newline
Hence, by taking the limit as $n\to \infty$, we conclude that (\ref{eq-conv}) holds whenever $f=P$ is a polynomial.
Now let us prove it for any $f\in C([a,b])$.

Corresponding to each $\epsilon >0$, by the Weierstrass Theorem,  there exists a polynomial $P^*$ such that
\[
\|f-P^*\|_\infty<\epsilon.
\]
Moreover, since $\lim_{n\rightarrow\infty}\|P^*- \tilde r(P^*)\|_\infty=0$, there exists $\nu_\epsilon >0$ such that
\[
\|P^*-\tilde r(P^*)\|_\infty <\epsilon, \qquad \forall n>\nu_\epsilon.
\]
Finally, note that by applying Thm. \ref{th-LC} with $h\sim h^*\sim n^{-1}$, we get
\[
\|\tilde r(F)\|_\infty \le  \Lambda_n\|F\|_\infty\le \C \|F\|_\infty, \qquad \forall F\in C([a,b]),
\]
where $\C>0$ is independent of $n$.

Summing up, from the above results we conclude that  $\forall \epsilon >0$ $\exists \nu_\epsilon >0$ such that $\forall n>\nu_\epsilon$ we have
\begin{eqnarray*}
\|f-\tilde r(f)\|_\infty &\le& \|f-P^*\|_\infty + \|P^*-\tilde r (P^*)\|_\infty +  \|\tilde r (f-P^*)\|_\infty\\
&\le & 2\epsilon + \Lambda_n \|f-P^*\|_\infty \le (2+\C) \epsilon
\end{eqnarray*}
that means (\ref{eq-conv}) holds for arbitrary $ f\in C([a,b])$.
\Proofend

In the following, we provide several estimates of the convergence order by supposing different degrees of smoothness for the function $f\in C([a,b])$.

Let us start by providing an error estimate in the case that $f$ is a H\"older continuous function satisfying
\begin{equation}\label{lip}
|f(x)-f(y)|\le M |x-y|^\alpha, \qquad \forall x,y\in [a,b],
\end{equation}
with $ 0<\alpha\le 1$ and $M>0$ independent of $x$ and $y$.

Denoted by $Lip_\alpha([a,b])$  the class of all such functions, we state the following
\begin{theorem}\label{th-lip}
Let $\tilde r(f,x)$ be the generalized FH interpolant corresponding to fixed parameters $d,\gamma\in\NN$, $\gamma \ge 2$,  arbitrary large $n\ge d$, and a distribution of nodes satisfying $h\sim h^*§\sim n^{-1}$.
For any $0<\alpha\le 1$, if we have $f\in Lip_\alpha([a,b])$ then we get
\begin{equation}\label{eq-lip}
|f(x)-\tilde r(f,x)|\le \C \left\{\begin{array}{lll}
\displaystyle\frac 1{n^\alpha} & \mbox{if} & \gamma >\alpha +1\\ [.19in]
\displaystyle\frac{\log n}{n^\alpha} & \mbox{if} & \gamma=\alpha +1
\end{array}\right.\qquad \forall x\in [a,b],
\end{equation}
where $\C>0$ is a constant as in Theorem \ref{th-LCrelated}.
\end{theorem}
{\it Proof of Theorem \ref{th-lip}}.
From (\ref{f-rf}) and (\ref{lip}) we deduce that
\[
|f(x)-\tilde r(f,x)|
 \le \sum_{k=0}^n \left|f(x)-f(x_k)\right| |b_k(x)|\le M  \sum_{k=0}^n \left|x-x_k\right|^\alpha |b_k(x)|
\]
holds  for all $x\in[a,b]$. Hence, the statement follows by applying Thm. \ref{th-LCrelated}.
\Proofend

Now, let us estimate the error in the case $f$ belongs to the class $C^{s}([a,b])$ of all functions that are $s$--times continuously differentiable in $[a,b]$.
\begin{theorem}\label{th-Cs}
Let $\tilde r(f,x)$ be the generalized FH interpolant corresponding to fixed parameters $d,\gamma\in\NN$ with $\gamma \ge 2$,  arbitrary large $n\ge d$, and a distribution of nodes satisfying $h\sim h^*§\sim n^{-1}$.
 For any integer $1\le s\le (d+1)$, if  $f\in C^s([a,b])$ then we have
\begin{equation}\label{eq-Cs}
|f(x)-\tilde r(f,x)|\le \C \left\{\begin{array}{lll}
\displaystyle\frac 1{n^s} & \mbox{if} & \gamma >s +1\\ [.19in]
\displaystyle\frac{\log n}{n^s} & \mbox{if} & \gamma=s +1 \\ [.19in]
\displaystyle\frac 1{n^{\gamma-1}} & \mbox{if} & 1<\gamma <s +1
\end{array}\right.\qquad \forall x\in [a,b],
\end{equation}
where $\C>0$ is a constant as in Theorem \ref{th-LCrelated}.
\end{theorem}
{\it Proof of Theorem \ref{th-Cs}}.
Due to the interpolation property (\ref{eq-interp}), it is sufficient to prove (\ref{eq-Cs}) in the case that
$x\in ]x_\ell, x_{\ell+1}[$ is arbitrarily fixed, for any $\ell\in\{0,\ldots, n-1\}$.

Since $f\in C^s([a,b])$ with $s\ge 1$, we can certainly consider the Taylor polynomial of $f$  centered at $x_\ell$ and having degree $s-1\ge 0$, namely
\begin{equation}\label{Taylor}
T(y)=\sum_{j=0}^{s-1}\frac{f^{(j)}(x_\ell)}{j!}(y-x_\ell)^j,\qquad \forall y\in [a,b].
\end{equation}
Recalling the Lagrange form of the remainder term, we get  the following error-bound
\begin{equation}\label{resto-Lag}
|f(y)-T(y)|\le \frac{\|f^{(s)}\|_\infty}{s!}\ |y-x_{\ell}|^{s} = \C  |y-x_{\ell}|^{s}  ,\qquad
\forall y\in [a,b]
\end{equation}
where here and in the following $\C>0$ denotes any constant as in Thm. \ref{th-LCrelated} which can take also different values at different occurrences.

Since $\deg T\le d$, we can use the polynomial preservation property (\ref{eq-inva}), that combined with (\ref{resto-Lag}) yields
\begin{eqnarray*}
 |f(x)-\tilde r(f,x)| &=&\left| f(x)-T(x) + \tilde r(T-f,x)\right|\\
&\le& |f(x)-T(x)|+\sum_{k=0}^n |T(x_k)-f(x_k)| |b_k(x)|\\
&\le& \C |x-x_\ell|^s + \C \sum_{k=0}^n |x_k-x_\ell|^s |b_k(x)|
\\ &\le&\frac \C{n^s}+\C \sum_{k=0}^n |x_k-x_\ell|^s |b_k(x)|
\end{eqnarray*}
where in the last inequality we used $|x-x_\ell|\le (x_{\ell+1}-x_\ell)\le h\sim n^{-1}$.\newline
On the other hand, since
\[
(u+v)^d\le 2^{d-1} \left(u^d+v^d\right),\qquad \forall u,v\in\RR^+,
\]
holds, we have
\[
|x_k-x_\ell|^s=|(x_k-x)+(x-x_\ell)|^s\le \left(|x_k-x|+|x-x_\ell|\right)^s
\le \C \left(|x_k-x|^s+|x-x_\ell|^s\right)
\]
Thus, we conclude that
\begin{eqnarray*}
 |f(x)-\tilde r(f,x)|
&\le& \frac\C{n^s}+ \C \sum_{k=0}^n |x_k-x_\ell|^s |b_k(x)|
\\ &\le&\frac\C{n^s}+\C \sum_{k=0}^n |x-x_k|^s |b_k(x)|+ \C |x-x_\ell|^s\sum_{k=0}^n |b_k(x)|\\
&\le& \frac\C{n^s} \left[1+\sum_{k=0}^n |b_k(x)|\right] +\C \sum_{k=0}^n |x-x_k|^s |b_k(x)|
\end{eqnarray*}
and the statement follows by applying Theorems \ref{th-LC} and \ref{th-LCrelated}.
\Proofend

We remark that Thm. \ref{th-Cs} for $s=d+1$ states that the analogous of (\ref{err-FH}) holds for generalized FH interpolants too, but under weaker assumption on the function $f$.

Now we estimate the error in the class $C^{s,\alpha}([a,b])$ of all functions  that are $s$--times continuously differentiable with $f^{(s)}\in Lip_\alpha([a,b])$, $0<\alpha\le 1$.
\begin{theorem}\label{th-Csa}
Let $\tilde r(f,x)$ be the generalized FH interpolant corresponding to fixed parameters $d,\gamma\in\NN$ with $\gamma \ge 2$,  arbitrary large $n\ge d$, and a distribution of nodes satisfying $h\sim h^*§\sim n^{-1}$. For any  $0<\alpha\le 1$ and each $s\in \NN$ with $s\le d$, if $ f\in C^{s,\alpha}([a,b])$ then  we have
\begin{equation}\label{eq-Csa}
|f(x)-\tilde r(f,x)|\le \C \left\{\begin{array}{lll}
\displaystyle\frac 1{n^{s+\alpha}} & \mbox{if} & \gamma >s +\alpha+1\\ [.19in]
\displaystyle\frac{\log n}{n^{s+\alpha}} & \mbox{if} & \gamma=s +\alpha+1 \\ [.19in]
\displaystyle\frac 1{n^{\gamma-1}} & \mbox{if} & 1<\gamma <s +\alpha+1
\end{array}\right.\qquad \forall x\in [a,b],
\end{equation}
where $\C>0$ is a constant as in Theorem \ref{th-LCrelated}.
\end{theorem}
{\it Proof of Theorem \ref{th-Csa}}.
The proof follows similarly to that of Thm. \ref{th-Cs}  but this time we take, as $T(y)$, the Taylor polynomial of $f$ having degree $s\le d$ and centered at $x_\ell$ (supposed that $x\in ]x_\ell,x_{\ell+1}[$). By (\ref{eq-inva}) we get
\[
|f(x)-\tilde r(f,x)|\le |f(x)-T(x)|+ \sum_{k=0}^n |f(x_k)-T(x_k)||b_k(x)|.
\]
Hence, taking into account Thm. \ref{th-LCrelated}, the crucial step to get the statement is proving that
\begin{equation}\label{T-err1}
|f(y)-T(y)|\le \C \ |y-x_{\ell}|^{s+\alpha} ,\qquad \forall y\in [a,b], \qquad \forall f\in C^{s,\alpha}([a,b])
\end{equation}
where $\C>0$ is independent of $y$ and $x_\ell$.

This can be easily proved by induction on $s$, keeping $0<\alpha\le 1$ arbitrarily fixed. \newline
Since the statement is trivial when $y=x_\ell$, let us arbitrarily fix,  for instance,  $x_\ell<y\le b$,  the proof is similar if $a\le y<x_\ell$.\newline
For $s=1$, (\ref{T-err1}) holds since by the Mean Value Theorem
\[
f(y)-f(x_\ell)=f^\prime(\xi)(y-x_\ell), \qquad x_\ell<\xi<y,
\]
and consequently $\forall f\in C^{1,\alpha}([a,b])$ we have
\begin{eqnarray*}
|f(y)-T(y)|&=&|f(y)-f(x_\ell)-f^\prime(x_\ell)(y-x_\ell)|\\ &=&|f^\prime(\xi)-f^\prime(x_\ell)||y-x_\ell|\\
&\le&\C|\xi-x_\ell|^{\alpha}|y-x_\ell|\le \C |y-x_{\ell}|^{\alpha+1}.
\end{eqnarray*}
Now suppose that (\ref{T-err1}) holds for $s-1$ and prove it for $s$.
By applying  the Cauchy Theorem to the following functions
\[
F(y):=f(y)-T(y), \qquad G(y):=(y-x_\ell)^{s+\alpha},
\]
we get there exists $\xi\in ]x_\ell,y[$ such that
\begin{eqnarray*}
\frac{|f(y)-T(y)|}{(y-x_\ell)^{s+\alpha}}&=&
\frac{|F(y)-F(x_\ell)|}{|G(y)-G(x_\ell)|}=\frac{|F^\prime(\xi)|}{|G^\prime(\xi)|}\\
&=&
\frac{\left|f^\prime(\xi)-
\left(\sum_{j=0}^{s-1}\frac{(f^\prime)^{(j)}(x_\ell)}{j!} (\xi-x_\ell)^j\right)\right|}
{(s+\alpha)|\xi-x_\ell|^{s+\alpha-1}}.
\end{eqnarray*}
Thus, since  (\ref{T-err1}) holds for $s-1$, by  applying it to $f^\prime\in C^{s-1,\alpha}([a,b])$, the previous estimate continues as follows
\[
\frac{|f(y)-T_(y)|}{(y-x_\ell)^{s+\alpha}}=
\frac{\left|f^\prime(\xi)-
\left(\sum_{j=0}^{s-1}\frac{(f^\prime)^{(j)}(x_\ell)}{j!} (\xi-x_\ell)^j\right)\right|}
{(s+\alpha)|\xi-x_\ell|^{s+\alpha-1}}\le \frac{\C |\xi-x_\ell|^{s-1+\alpha}}{(s+\alpha)|\xi-x_\ell|^{s+\alpha-1}}\le \C
\]
 i.e., (\ref{T-err1}) holds for $s$ too.
\Proofend

Finally, we focus on the general case that $f:[a,b]\to \RR$ is a function of bounded variation ($f\in BV([a,b])$) and show that all the previous error estimates continue to hold locally, in all compact subintervals $I$ where we have the above-prescribed smoothness, namely $f\in Lip_\alpha(I)$ with $0<\alpha\le 1$, or $f\in C^s(I)$ with $s\le d+1$, or $f\in C^{s,\alpha}(I)$, with $s\le d$. In order to give a unified treatment of all these cases, we introduce the following notation
\begin{equation}\label{hp-BV}
C^{s,\alpha}:=\left\{\begin{array}{lll}
Lip_\alpha & \mbox{if $0<\alpha\le 1$} & \mbox{and $s=0$}\\
C^s & \mbox{if $\alpha=0$} & \mbox{and $s\in\NN$ with $s\le d+1$}\\
C^{s,\alpha} & \mbox{if $0<\alpha\le 1$} & \mbox{and $s\in\NN$ with $s\le d$}
\end{array}\right.
\end{equation}
\begin{theorem}\label{th-local}
Let $f\in BV([a,b])$ and let $\tilde r(f,x)$ be the generalized FH interpolant of $f$ corresponding to fixed $d,\gamma\in\NN$ with $\gamma\ge 2$, arbitrarily large $n\ge d$ and  nodes satysfying $h\sim h^*\sim n^{-1}$. Under the assumptions in (\ref{hp-BV}), if $f\in C^{s,\alpha}(I)$  with $I=[a^\prime, b^\prime]\subset [a,b]$ , then we have
\begin{equation}\label{eq-local}
|f(x)-\tilde r(f,x)|\le \C\left\{
\begin{array}{lll}
\displaystyle\frac 1{n^{s+\alpha}} & \mbox{if} & \gamma >s +\alpha+1\\ [.19in]
\displaystyle\frac{\log n}{n^{s+\alpha}} & \mbox{if} & \gamma=s +\alpha+1 \\ [.19in]
\displaystyle\frac 1{n^{\gamma-1}} & \mbox{if} & 1<\gamma <s +\alpha +1
\end{array}\right.\qquad \forall x\in  [a^\prime, b^\prime]
\end{equation}
where $\C>0$ is a constant as in Theorem \ref{th-LCrelated}.
\end{theorem}
{\it Proof of Theorem \ref{th-local}}.
As usual, we consider the case $x\in [a^\prime, b^\prime]-\{x_0,\ldots, x_n\}$ is arbitrarily fixed and we denote by $\C$ a positive constant that can take different values at different occurrences and has the qualities described in Thm.\ref{th-LCrelated}.\newline As $n\to\infty$, we can suppose that $x_\ell<x<x_{\ell+1}$ with $[x_\ell,x_{\ell+1}]\subset [a^\prime, b^\prime]$. In this way, we can take the Taylor polynomial of $f$ centered at $x_\ell$ and with degree
\[\left\{ \begin{array}{ll}
s &\mbox{if $s\le d$}\\
s-1 & \mbox{if $s=d+1$ and $\alpha=0$}
\end{array}\right.\]
Denoted by $T(y)$ such polynomial, since $f\in C^{s,\alpha}([a^\prime,b^\prime])$, we have (cf. (\ref{resto-Lag}) and (\ref{T-err1}) )
\begin{equation}\label{T-err}
|f(y)-T(y)|\le \C \ |y-x_{\ell}|^{s+\alpha} ,\qquad \forall y\in [a^\prime,b^\prime],
\end{equation}
and since $\deg(T)\le d$, by (\ref{eq-inva}), we also have $\tilde r(T)=T$.
\newline
Consequently, we get
\begin{eqnarray*}
|f(x)-\tilde r(f,x)|&\le&|f(x)-T(x)|+|\tilde r(T-f,x)|\\
&\le& \C |x-x_\ell|^{s+\alpha}+ \sum_{k=0}^n |T(x_k)-f(x_k)| |b_k(x)|,
\end{eqnarray*}
and  setting
\begin{eqnarray*}
I_1&=&\{k\in\{0,1,\ldots, n\} :\ |x-x_k|\le (b^\prime- a^\prime)\}, \\
I_2&=&\{k\in\{0,1,\ldots, n\} :\ |x-x_k|> (b^\prime- a^\prime)\}
\end{eqnarray*}
we obtain
\begin{eqnarray*}
|f(x)-\tilde r(f,x)|&\le&\frac\C{n^{s+\alpha}}+ \left\{\sum_{k\in I_1} + \sum_{k\in I_2} \right\} |T(x_k)-f(x_k)| |b_k(x)|\\
&=:&\frac\C{n^{s+\alpha}} + s_1(x)+s_2(x)
\end{eqnarray*}
As regards $s_1(x)$, the estimate (\ref{T-err}) can be applied because $k\in I_1$ implies $x_k\in [a^\prime, b^\prime]$. Consequently, by  (\ref{T-err}), Thm. \ref{th-LC} and Thm.\ref{th-LCrelated}, we get
\begin{eqnarray*}
s_1(x)&:=&\sum_{k\in I_1}  |T(x_k)-f(x_k)| |b_k(x)|\\
&\le& \C \sum_{k\in I_1}  |x_k-x_\ell|^{s+\alpha} |b_k(x)|\\
&\le& \C \sum_{k=0}^n \left(|x-x_\ell|^{s+\alpha}+ |x-x_k|^{s+\alpha}\right) |b_k(x)|\\
&\le& \frac{\C}{n^{s+\alpha}} +    \C \sum_{k=0}^n|x-x_k|^{s+\alpha}|b_k(x)|
\le \C\left\{
\begin{array}{lll}
\displaystyle\frac 1{n^{s+\alpha}} & \mbox{if} & \gamma >s +\alpha+1\\ [.1in]
\displaystyle\frac{\log n}{n^{s+\alpha}} & \mbox{if} & \gamma=s +\alpha+1 \\ [.1in]
\displaystyle\frac 1{n^{\gamma-1}} & \mbox{if} & 1<\gamma <s +\alpha +1
\end{array}\right.
\end{eqnarray*}
Finally, concerning $s_2(x)$, we note that
\begin{eqnarray*}
s_2(x)&:=& \sum_{k\in I_2} \left[\frac{|T(x_k)-f(x_k)|}{|x-x_k|^{s+\alpha}}\right] |x-x_k|^{s+\alpha} |b_k(x)|\\
&\le& \left[\sup_{k\in I_2}\frac{|T(x_k)-f(x_k)|}{|x-x_k|^{s+\alpha}}\right]
 \sum_{k\in I_2} |x-x_k|^{s+\alpha} |b_k(x)|\\
&\le& \frac{\|T-f\|_\infty}{(b'-a') ^{s+\alpha}}\sum_{k\in I_2} |x-x_k|^{s+\alpha} |b_k(x)|\\
&\le& \C \sum_{k=0}^n|x-x_k|^{s+\alpha}|b_k(x)|
\end{eqnarray*}
Hence,  by Thm. \ref{th-LCrelated}, we obtain that $s_2(x)$ can be estimate as $s_1(x)$.
\Proofend
\section{Some numerical experiments}\label{sec:ex}
In this section we compare the error function
$$
e_{n,d,\gamma}(x) = f(x) -\tilde r(f,x)
$$
for different functions $f$ and different values of $n$, $d$ and $\gamma$.
We denote the maximum error by
$$
E_{n,d,\gamma} = \max_{x \in [-1,+1]} | e_{n,d,\gamma} (x)|.
$$
As interpolation nodes, we take equidistant points in the interval $[-1,+1]$.
Note that the generalized FH interpolants are defined for a general configuration of the interpolation points.

\begin{figure}[!htb]%
\begin{center}
\includegraphics[scale = 0.6]{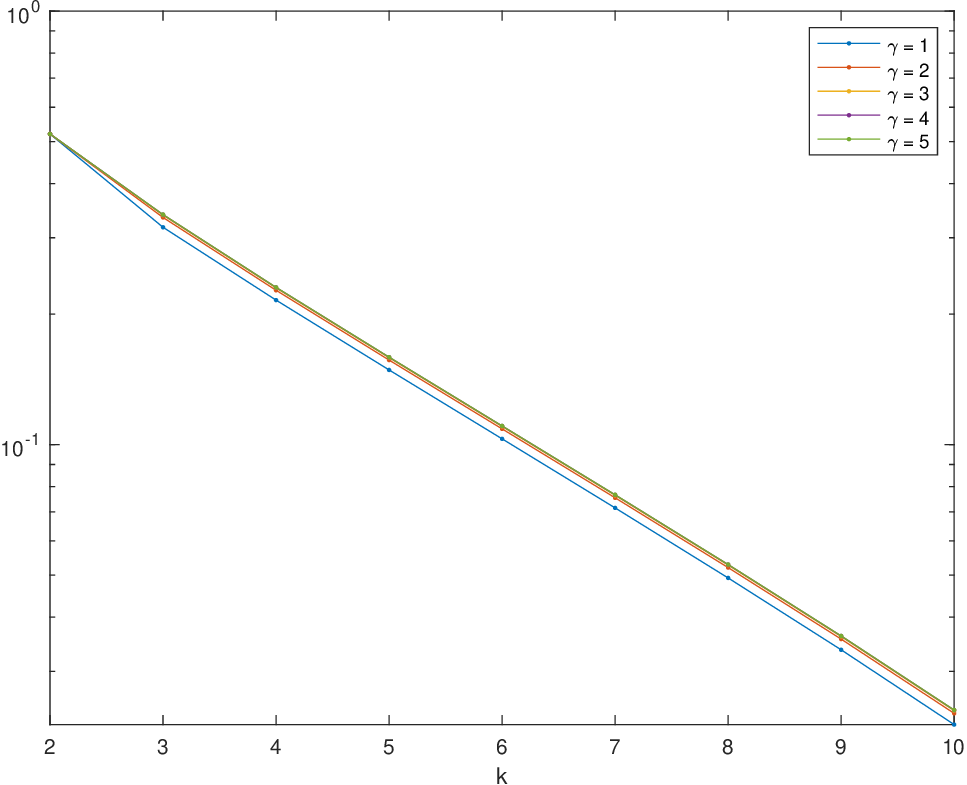}
\end{center}
\caption{ Maximum error $E_{2^k,d,\gamma}$ for $k=1,2,\ldots,10$ and $\gamma = 1,\ldots,5$ with $d = 2$ and $f(x) = |x|^{0.5}$.}\label{fig101}
\end{figure}

\begin{figure}[!htb]%
\begin{center}
\includegraphics[scale = 0.6]{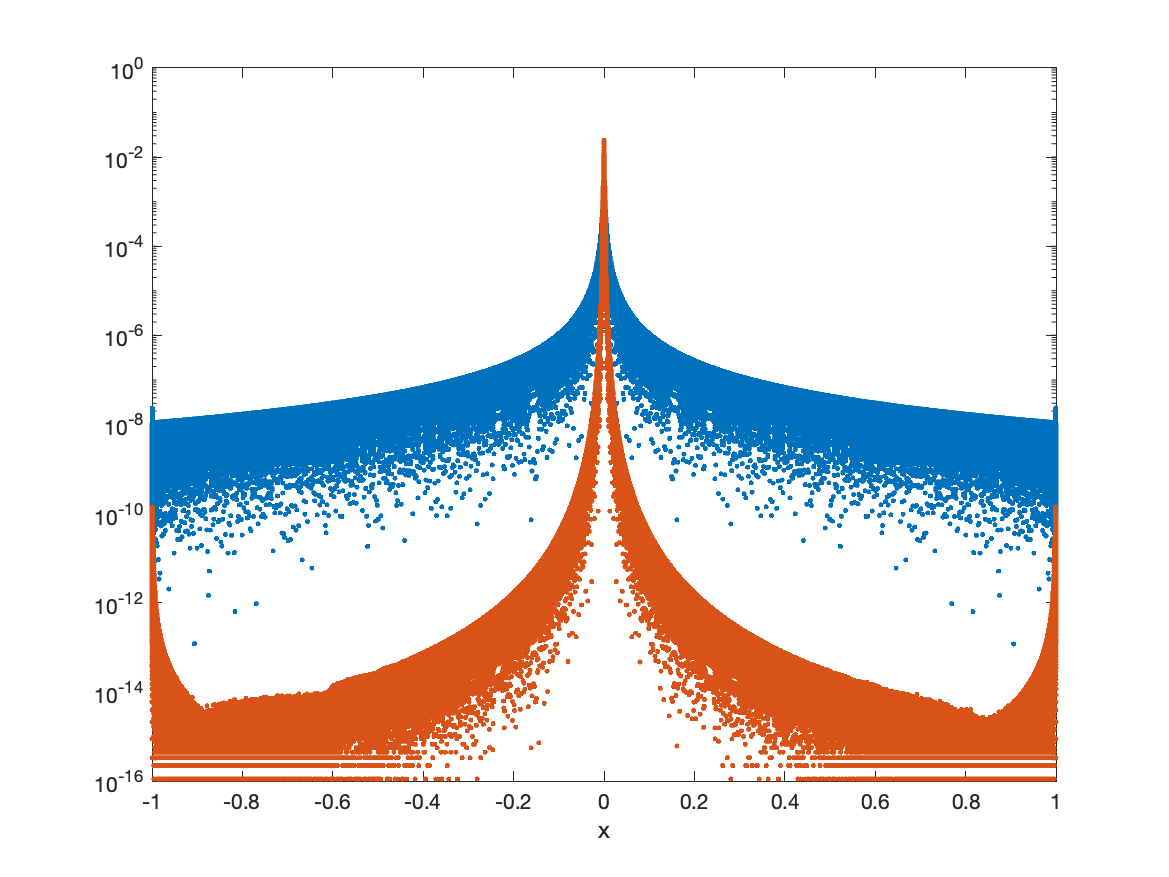}
\end{center}
\caption{ The error function $e_{n,d,\gamma}$ for $n = 1024$, $d=2$ and $\gamma = 1,2$, for the function $f(x) = |x|^{0.5}$.
The blue dots represent the error function for $\gamma = 1$ (original FH interpolants), the red dots for $\gamma =2$. }\label{fig104}
\end{figure}

{\bf Experiment 1:} Let us consider the function $f(x) = |x|^{0.5} \in Lip_\alpha ([-1,+1])$ with $\alpha = 0.5$ to illustrate
Thm.~\ref{th-lip}. According to this theorem,  $\gamma$ should be taken greater than $\alpha+1 = 1.5$.
The numerical experiment shows that in this case the theorem is also valid when $\gamma = 1$.
We get that for $n = 2^k$, $k=1,2,\ldots,10$,  the maximum error $E_{2^k,d,\gamma}$ is divided by a factor approximately equal to $\sqrt{2}$
when $k$ is increased by one.  So the maximum error behaves as $\O(h^{0.5})$.
This is true for $\gamma = 1,2,\ldots$ and $d = 0,1,2,\ldots$.
Figure \ref{fig101} illustrates this by plotting $E_{2^k,d,\gamma}$ for $k=1,2,\ldots,10$ and $\gamma = 1,\ldots,5$ with $d = 2$.
The curve for $\gamma$ is above the one for $\gamma-1$.
So, one would think that increasing the value of $\gamma$ is not a good idea.
However, the factors between the errors for $\gamma=5$ and $\gamma=1$ are small and around $1.08$.
More importantly the error function $e_{n,d,\gamma}$ for this function $f(x) = |x|^{0.5}$
behaves much better for  $\gamma=2$ compared to $\gamma=1$.
Figure \ref{fig104} shows the error function $e_{n,d,\gamma}$ for $n = 1024$, $d=2$ and $\gamma = 1,2$.
The blue dots represent the error function for $\gamma = 1$ (original FH interpolants), the red dots for $\gamma =2$.
Although the peak increases slightly for increasing values of $\gamma$, the behaviour away from the origin is
much better.  We did not show the behaviour of the error function for $\gamma = 3,4,5$.
In this case the behaviour is between the behaviour for $\gamma = 1$ and
$\gamma=2$.

\begin{figure}[!htb]%
\begin{center}
\includegraphics[scale = 0.6]{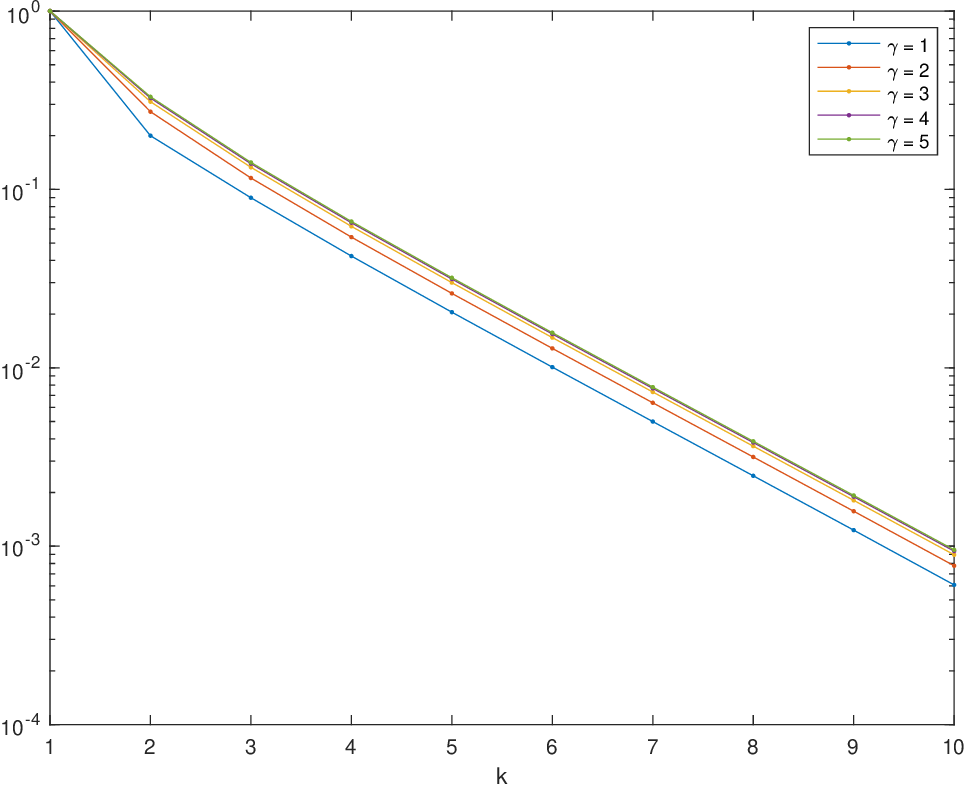}
\end{center}
\caption{ Maximum error $E_{2^k,d,\gamma}$ for $k=1,2,\ldots,10$ and $\gamma = 1,\ldots,5$ with $d = 1$ and $f(x) = |x|$.}\label{fig01}
\end{figure}

\begin{figure}[!htb]%
\begin{center}
\includegraphics[scale = 0.6]{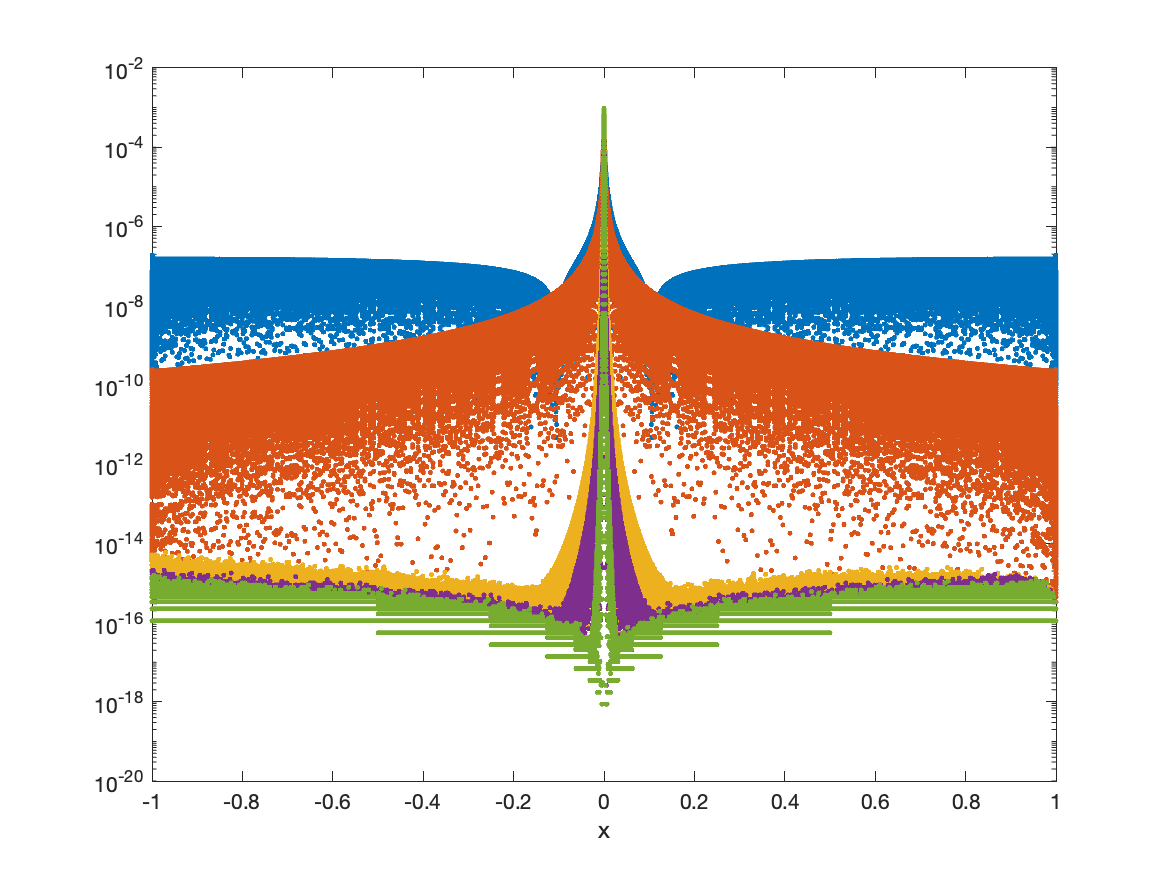}
\end{center}
\caption{ The error function $e_{n,d,\gamma}$ for $n = 1024$, $d=1$ and $\gamma = 1,2,\ldots,5$, for the function $f(x) = |x|$.
The blue dots represent the error function for $\gamma = 1$ (original FH interpolants), the red dots for $\gamma =2$
and so on. }\label{fig02}
\end{figure}

{\bf Experiment 2:} Let us consider the function $f(x) = |x| \in Lip_\alpha ([-1,+1])$ with $\alpha = 1$.
The conditions of Thm. \ref{th-lip} require that $\gamma$ should be greater than $\alpha+1 = 2$ to obtain a maximum error
that behaves as $\O(h)$.
However, the numerical experiments indicate that this behaviour  is also obtained for $\gamma =1$ and $2$.
We have that for $n = 2^k$, $k=1,2,\ldots,10$,  the maximum error $E_{2^k,d,\gamma}$ is divided by a factor approximately equal to $2$
when $k$ is increased by one.  So the maximum error behaves as $\O(h)$.
This is true for $\gamma = 1,2,\ldots$ and $d = 0,1,2,\ldots$.
Figure \ref{fig01} illustrates this by plotting $E_{2^k,d,\gamma}$ for $k=1,2,\ldots,10$ and $\gamma = 1,\ldots,5$ with $d = 1$.
The curve for $\gamma$ is above the one for $\gamma-1$.
Also here, one would think that increasing the value of $\gamma$ is not a good idea.
However, the factors between the errors for $\gamma=5$ and $\gamma=1$ are also small and around $1.55$.
More importantly,  the error functions $e_{n,d,\gamma}$  in the case $f(x) = |x|$
also behave much better for increasing values of $\gamma$.
Figure \ref{fig02} shows the error function $e_{n,d,\gamma}$ for $n = 1024$, $d=1$ and $\gamma = 1,2,\ldots,5$.
The blue dots represent the error function for $\gamma = 1$ (original FH interpolants), the red dots for $\gamma =2$
and so on.
Although the peak increases slightly for increasing values of $\gamma$, the behavior away from the origin is
much better.

\begin{figure}[!htb]%
\begin{center}
\includegraphics[scale = 0.6]{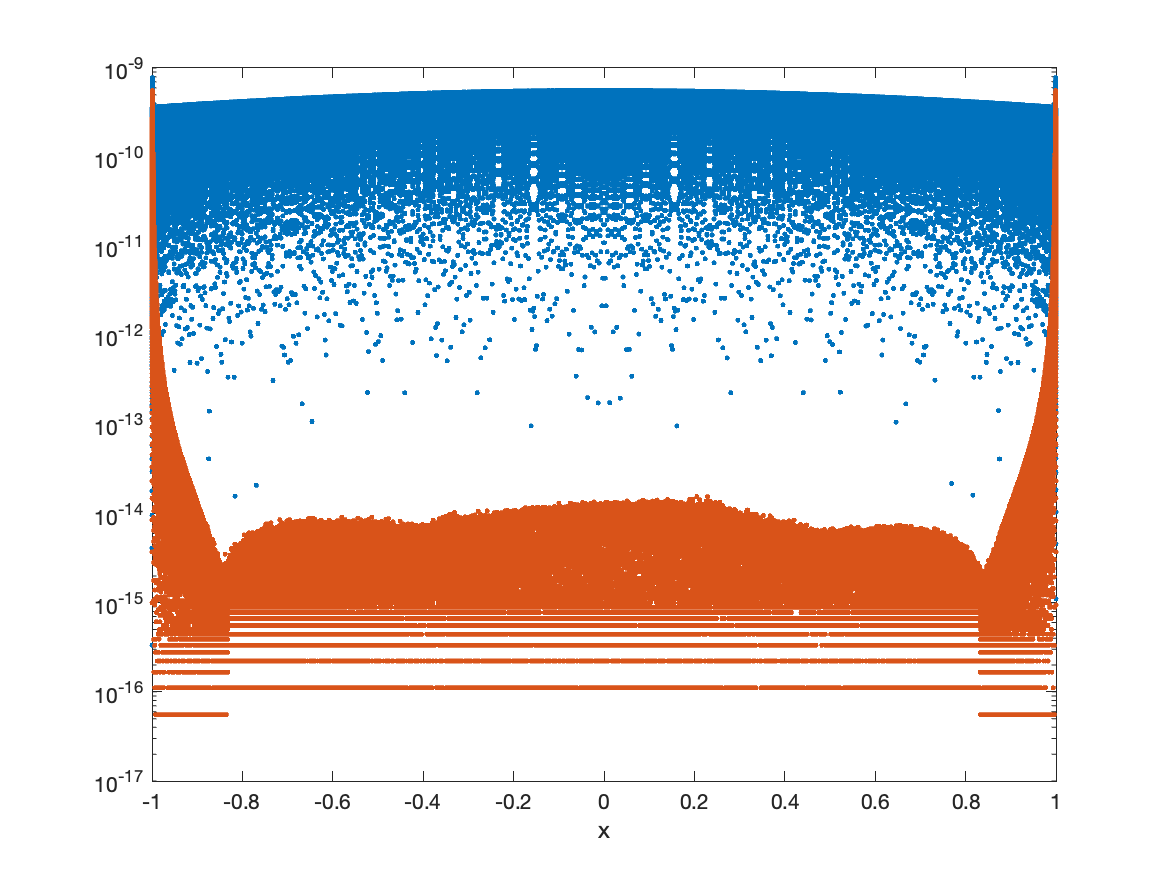}
\end{center}
\caption{ The error function $e_{n,d,\gamma}$ for $n = 1024$, $d=2$ and $\gamma = 1,2$, for the function $f(x) = e^{-x^2}$.
The blue dots represent the error function for $\gamma = 1$ (original FH interpolants), the red dots for $\gamma =2$.
 }\label{fig04}
\end{figure}

{\bf Experiment 3:} Let us consider the analytic function $f(x) = e^{-x^2} \in \C^{\infty}([-1,+1])$.
Figure \ref{fig04} shows the error function $e_{n,d,\gamma}$ for $n = 1024$, $d=2$ and $\gamma = 1,2$.
The blue dots represent the error function for $\gamma = 1$ (original FH interpolants), the red dots for $\gamma =2$.
{\color{black} Following the theory, the maximum error for the interpolants should behave as $\O(h^{d+1})$ when $\gamma=1$ (cf. (\ref{err-FH})). This can be also observed in practice for $\gamma=1$ but for $\gamma=2$ too, although Thm.~\ref{th-Cs} predicts a lower order of convergence when $1<\gamma\le d+2$.}
The behavior in between the peaks is much better when $\gamma =2$ compared to $\gamma=1$.

\begin{figure}[!htb]%
\begin{center}
\includegraphics[scale = 0.6]{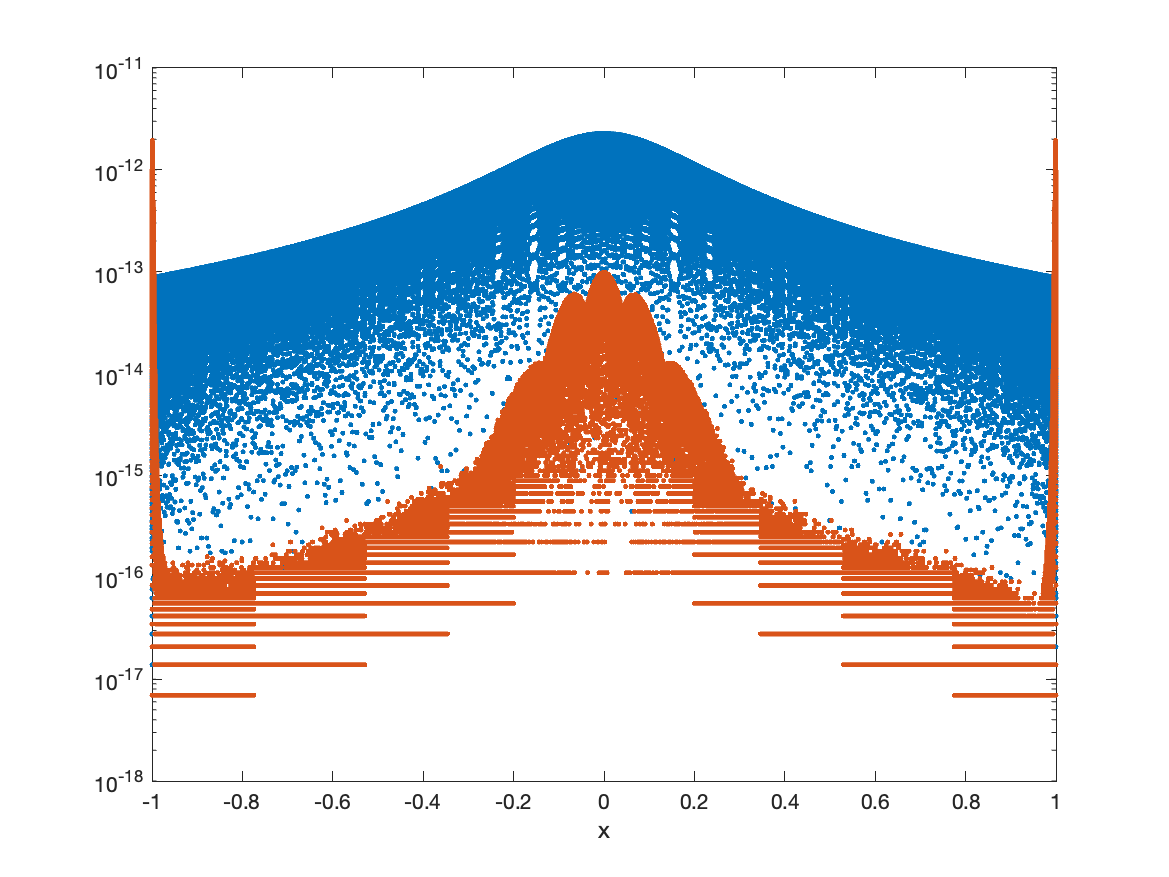}
\end{center}
\caption{ The error function $e_{n,d,\gamma}$ for $n = 1024$, $d=2$ and $\gamma = 1,2$, for the function $f(x) = 1/(1+25x^2)$.
The blue dots represent the error function for $\gamma = 1$ (original FH interpolants), the red dots for $\gamma =2$.
 }\label{fig206}
\end{figure}

\begin{figure}[!htb]%
\begin{center}
\includegraphics[scale = 0.6]{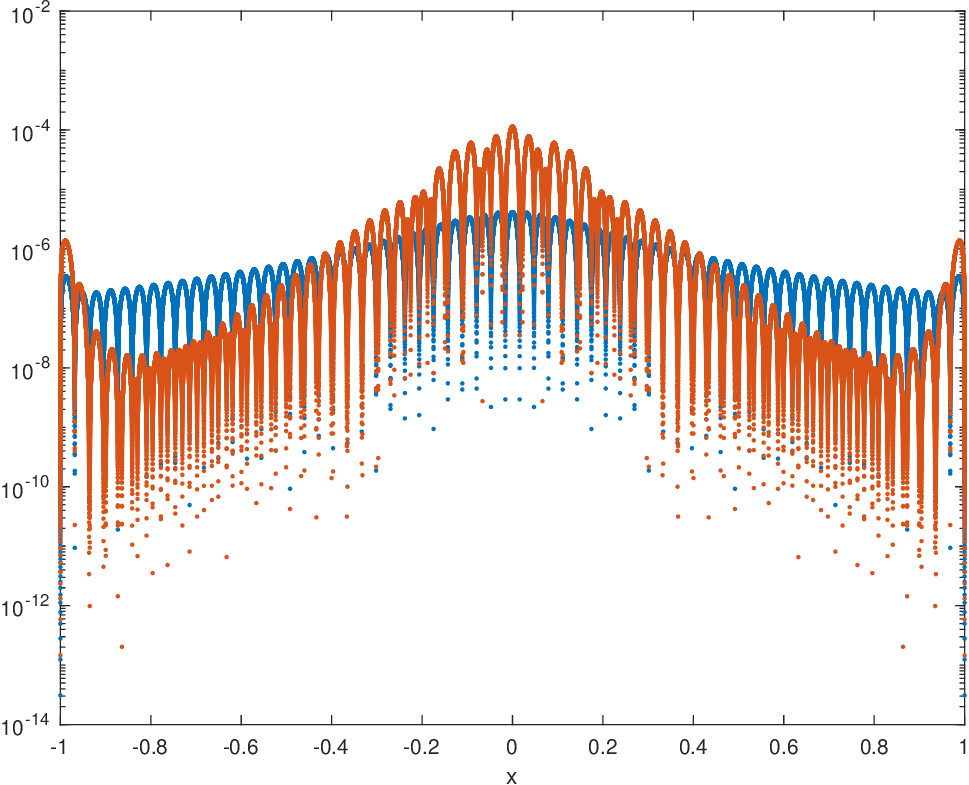}
\end{center}
\caption{ The error function $e_{n,d,\gamma}$ for $n = 64$, $d=2$ and $\gamma = 1,2$, for the function $f(x) = 1/(1+25x^2)$.
The blue dots represent the error function for $\gamma = 1$ (original FH interpolants), the red dots for $\gamma =2$.
 }\label{fig07}
\end{figure}

{\bf Experiment 4:} Let us consider the analytic Runge-function $f(x) = 1/(1+25x^2) \in \C^{\infty}([-1,+1])$.
Figure \ref{fig206} shows the error function $e_{n,d,\gamma}$ for $n = 1024$, $d=2$ and $\gamma = 1,2$.
The blue dots represent the error function for $\gamma = 1$ (original FH interpolants), the red dots for $\gamma =2$.
{\color{black} As in the previous experiment, according to (\ref{err-FH}), the maximum error for the interpolants behaves as $\O(h^{d+1})$ when $\gamma=1$, and this is also observed for $\gamma=2$.}
Moreover, the behavior in between the peaks is much better when $\gamma =2$ compared to $\gamma=1$.
However, if we take a similar figure for $n = 2^6 = 64$ we obtain Figure~\ref{fig07}.
In this case, $\gamma=2$ performs worse compared to $\gamma=1$.

\begin{figure}[!htb]%
\begin{center}
\includegraphics[scale = 0.4]{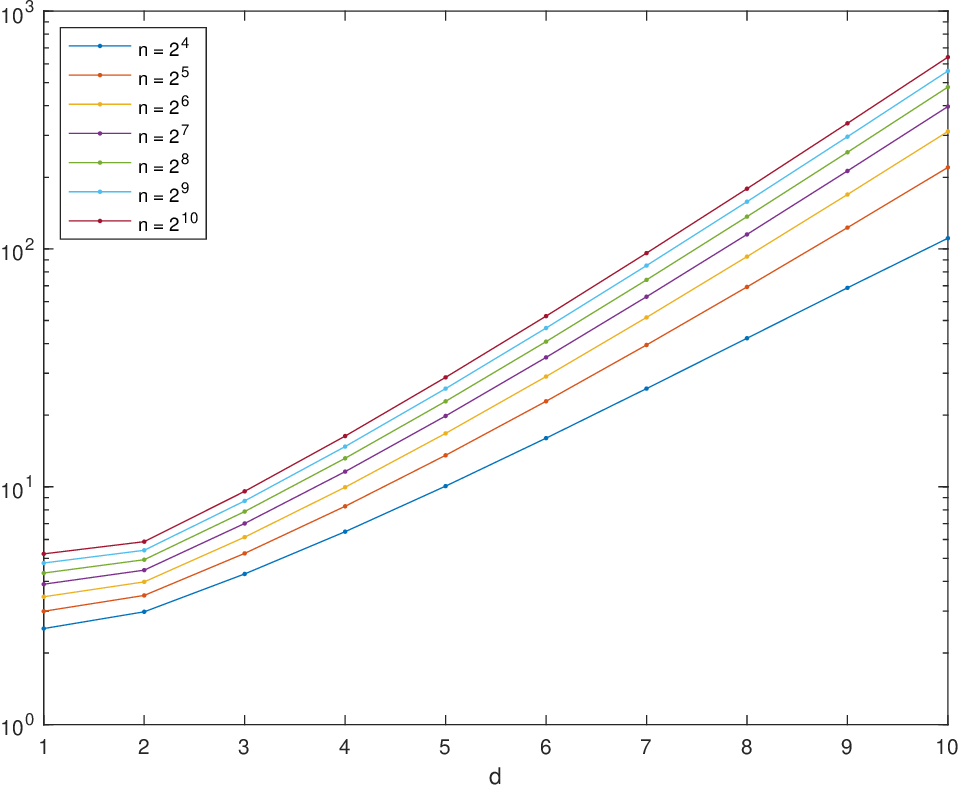}
\includegraphics[scale = 0.4]{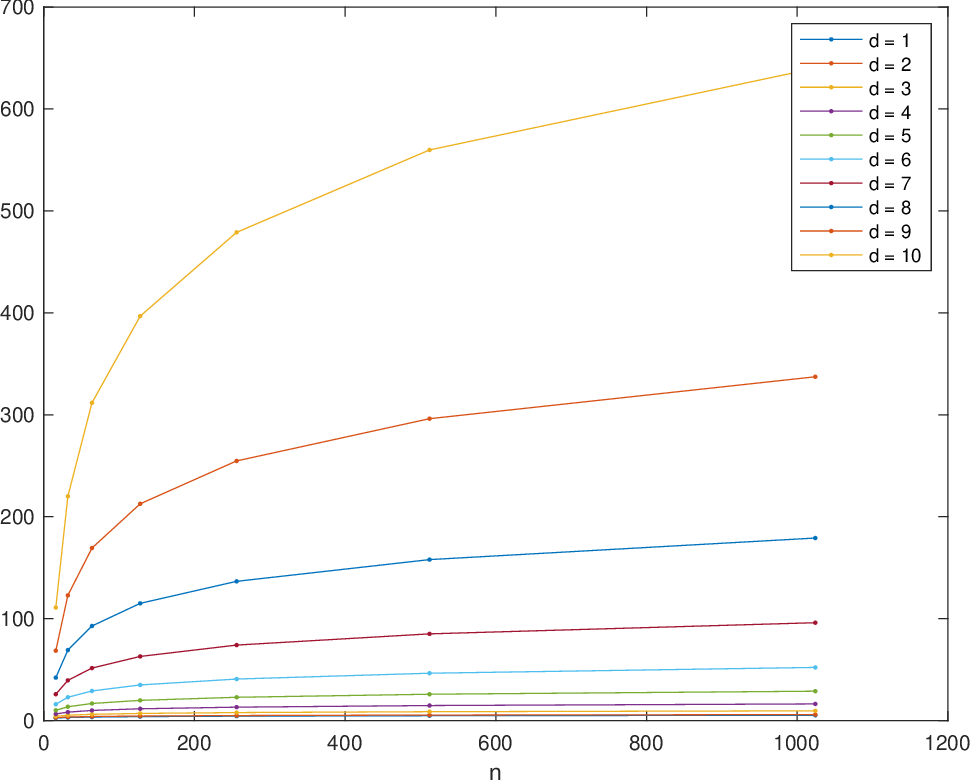}
\end{center}
\caption{ The Lebesgue constant for $\gamma = 1$ as function of $d$ for $n=2^k,k=4,5,\ldots,10$ (left) and as function of $n$ for $d=1,2,\ldots,10$ (right).
 }\label{fig0809}
\end{figure}

\begin{figure}[!htb]%
\begin{center}
\includegraphics[scale = 0.4]{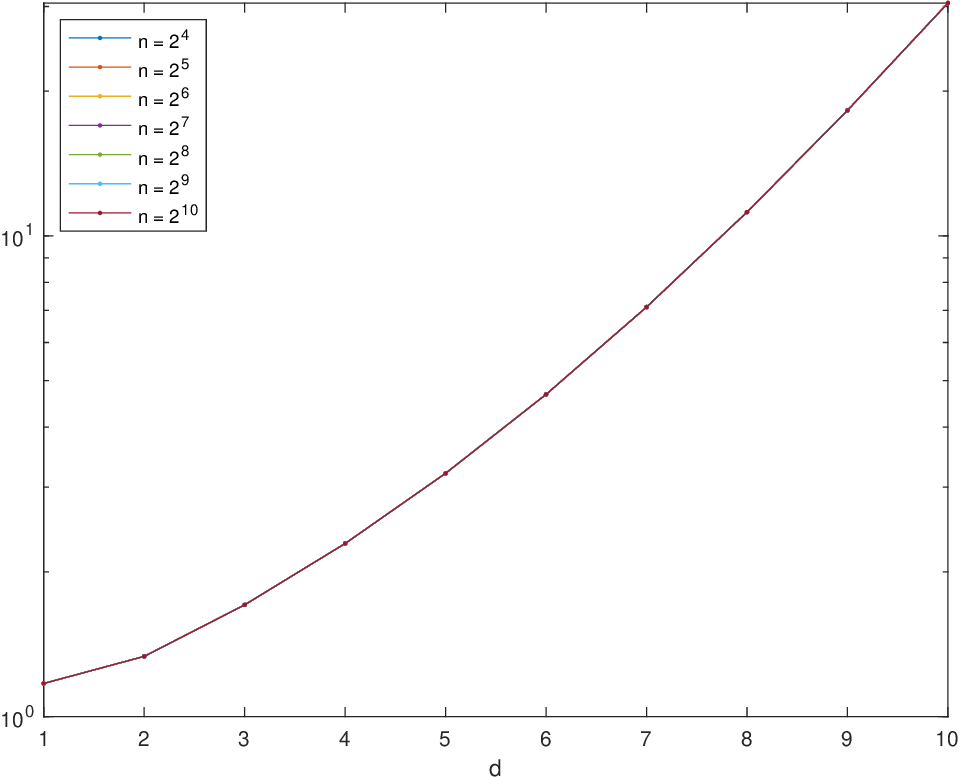}
\includegraphics[scale = 0.4]{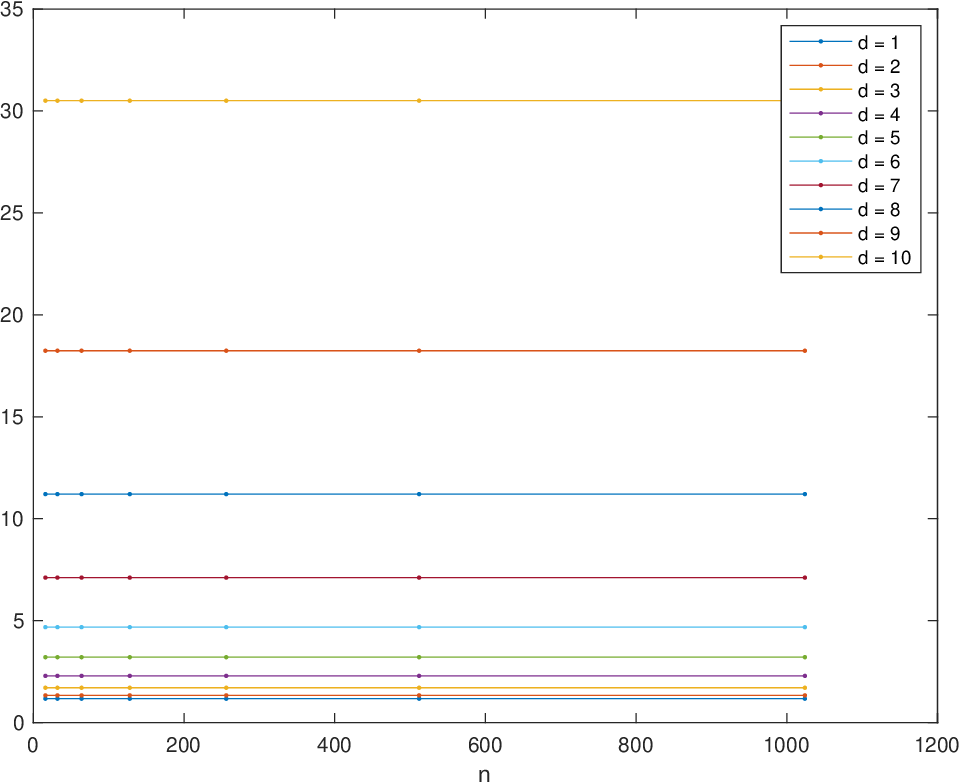}
\end{center}
\caption{ The Lebesgue constant for $\gamma = 2$ as function of $d$ for $n=2^k,k=4,5,\ldots,10$ (left) and in function of $n$ for $d=1,2,\ldots,10$ (right).
 }\label{fig1011}
\end{figure}

\begin{figure}[!htb]%
\begin{center}
\includegraphics[scale = 0.6]{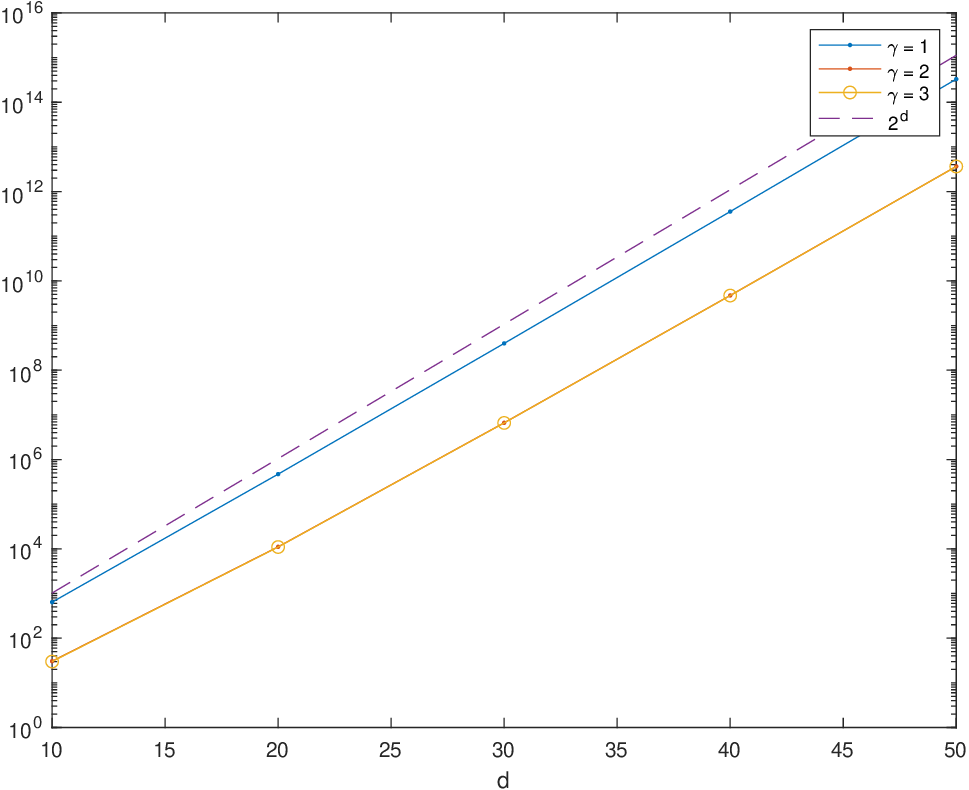}
\end{center}
\caption{ The Lebesgue constant for $d=10:10:50$, $n=2^k, k=10$ and for $\gamma=1,2,3$.
We compare this behaviour with the plot $2^d$.
 }\label{fig14}
\end{figure}

{\bf Experiment 5:} In \cite{hormann2012barycentric} an upper bound is derived for the Lebesgue constant in case $\gamma=1$:
\begin{equation}\label{eq:lebcon}
\Lambda_n \leq 2^{d} (1+\log n).
\end{equation}
To illustrate this, in Figure~\ref{fig0809} the Lebesgue constant is plotted in function of $d$ where the different lines correspond to values of $n = 2^k, k=4,5,\ldots,10$
(left) and in function of $n$ for different values of $d=1,2,\ldots,10$ (right).
The left figure clearly demonstrates the factor $2^{d}$ in (\ref{eq:lebcon}) while the right figure illustrates the $\log n$ dependency.
Plotting similar figures for $\gamma=2$ we obtain Figure~\ref{fig1011}.
This figure shows that for $\gamma=2$ the Lebesgue constant is independent of $n$.
For $\gamma=3$ we obtain similar figures.
In Figure~\ref{fig14} we plot the Lebesgue constant for $d=10:10:50$, $n=2^k, k=10$ and for $\gamma=1,2,3$.
We compare this behaviour with the plot $2^d$.
This figure shows that the Lebesgue constant behaves as $\C 2^d$ with $\C$ independent of $d$.

\begin{figure}[!htb]%
\begin{center}
\includegraphics[scale = 0.6]{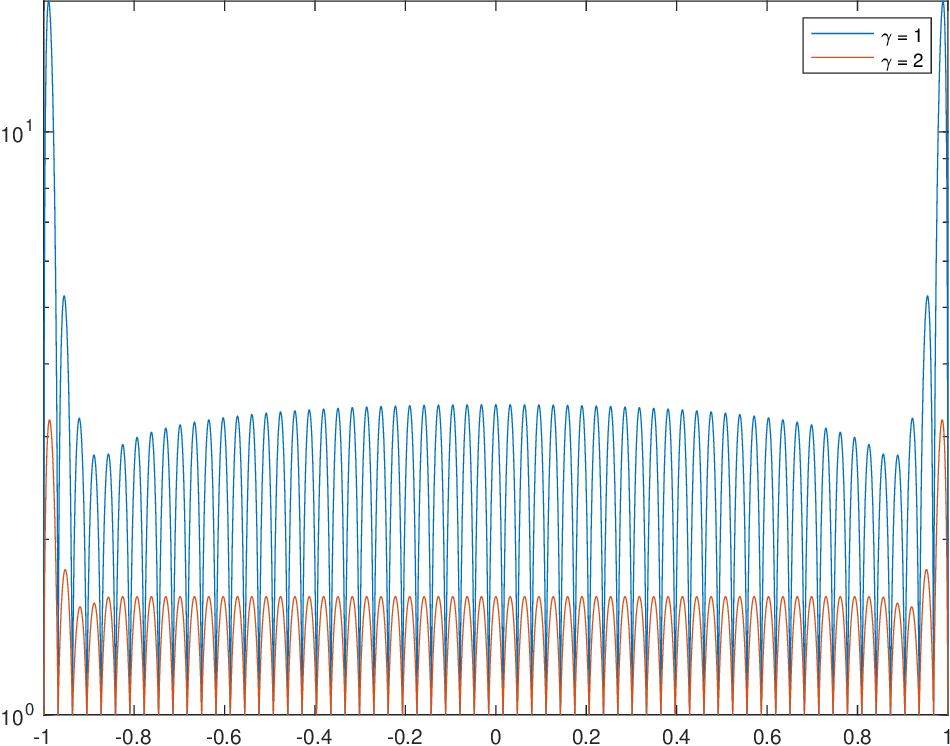}
\end{center}
\caption{ The Lebesgue function for $n = 2^6$, $d=5$ and $\gamma = 1$ and $2$.
 }\label{fig301}
\end{figure}

In Figure~\ref{fig301} the Lebesgue function is plotted for $n = 2^6$, $d = 5$ and $\gamma =1$ and $2$.

{\bf Experiment 6:} In this experiment, we compare the elapsed time of the classical Floater-Hormann algorithm
with computing and evaluating the new approximant.
We implemented the two algorithms in Matlab and measured the elapsed time using the {\tt tic} - {\tt toc} commands.
We used the following values: $d = 5$,  $n = 2^{10}$, $\gamma = 3$ and $m = 10^5$, the number of $x$-values in which the approximant
is evaluated.
The weights for the classical barycentric form can be computed beforehand by $\O(n d^2)$ FLOPS using
(\ref{bar-classic}) or in $\O(n d)$ FLOPS using a more complicated pyramid algorithm \cite{hormann2016pyramid}.
We implemented the former method.
Running the methods $50$ times and averaging, computing the weights for the classical FH approximant took $4\cdot  10^{-7} n d^2$ seconds while evaluating it in the $m$ points
costed $2\cdot 10^{-9} m n $ seconds.
Computing and evaluating the new approximant took $2 \cdot 10^{-9} m n d^2$ seconds.
For a different value of $\gamma$, we obtained comparable results.

\section{Conclusion}\label{sec:con}
In this paper, we have defined a whole family of generalized Floater--Hormann interpolants depending on an additional parameter $\gamma \in \NN$, besides the usual parameter $d\in\NN$.
For $\gamma=1$ we obtain the original Floater--Hormann interpolants.
The numerical examples show that this family has potential to approximate non-smooth as well as smooth functions.
In future work, the (sub-)optimal choice of the parameters $d$ and $\gamma$ could be investigated when the $n$ interpolation points are given.
For the original Floater--Hormann interpolants G\"uttel and Klein \cite{guttel2012convergence} have developed a heuristic method to determine the parameter $d$.
To remedy the bad behaviour of the error function at the endpoints, Klein \cite{klein2013extension} designed a method adding some interpolation points at the two endpoints.
A similar technique could be applied to the generalized Floater--Hormann approximants here introduced.
However, this approach seems useful only when the interpolated function is periodic \cite{camargo2017comparison}.
Moreover, in \cite{camargo2017stability} it is shown that Klein's method is numerically unstable. It could be interesting to investigate the well-conditioned case $d=0$ for all $\gamma\ge 2$, looking for some improvement w.r.t. the cases $\gamma=1,2$ studied in \cite{berrut1988rational, zhang2022rational}.
This and the following questions are left for future research: (A) Refine the theoretical bound of the Lebesgue constant in (\ref{eq-LC}) and also state a lower bound according to the numerical experiments; (B) Investigate the Lebesgue constants and the error for other configurations of nodes too; (C) Deeper explore the role of $\gamma$ in view of the numerical experiments that suggest larger theoretical bounds on $\gamma$.

\section*{Acknowledgments}
{The authors wish to thank the anonymous reviewers for their valuable comments.}

\end{document}